\documentclass[12pt,reqno]{amsart}
\usepackage[T1]{fontenc}
\usepackage{lmodern,microtype,graphicx}
\usepackage[margin = 1.2in]{geometry}
\usepackage{amsmath,amssymb,amsthm,enumitem}
\usepackage{booktabs,longtable,array,float,caption,tikz}
\usepackage[colorlinks = true,citecolor = black,linkcolor = black,urlcolor = black,hypertexnames = false]{hyperref}
\numberwithin{equation}{section}
\newtheorem{theorem}{Theorem}
\newtheorem{lemma}[theorem]{Lemma}
\newtheorem{proposition}[theorem]{Proposition}
\newtheorem{corollary}[theorem]{Corollary}
\newtheorem{problem}{Problem}

\newcommand{\R}{\mathbb R}
\newcommand{\F}{\mathbb F}
\newcommand{\EE}{\mathbb E_\mu}
\newcommand{\SO}{\mathrm{SO}}

\newcommand{\one}{\mathbf 1}

\title[Euclidean Ramsey Theory for Lines and Finite Configurations]{Euclidean Ramsey Theory for Lines and Finite Configurations}
\author{Ryan Hoff}
\address{Oberlin College}
\email{hoffcryan@gmail.com}
\author{Gary Hu}
\address{Williams College}
\email{gh7@williams.edu}
\author{Yumo Liu}
\address{Bard College}
\email{yy2323@bard.edu}
\date{}

\begin{document}

\begin{abstract}
  A classical question in Euclidean Ramsey theory, posed by Erdős, Graham, Montgomery, Rothschild, Spencer, and Straus, asks which finite configurations $X\subset\mathbb E^2$ must appear in blue under every red--blue coloring of the plane having no red pair of points distance one apart. For $m \ge 2$, let $\ell_m$ denote the line of $m$ points with consecutive points distance one apart. Tsaturian showed that $\ell_5$ is unavoidable; we improve this to $\ell_7$. Additionally, we show that in every three-coloring of the plane in which two colors contain no $\ell_2$, the third color must contain an $\ell_3$. To conclude our discussion on lines, we prove several related results on two-colorings of higher-dimensional Euclidean spaces where both colors avoid longer lines. Our methods also extend to other point configurations. In any two-coloring of the plane in which one color contains no $\ell_2$, we show that the other color must contain both a regular pentagon of side length one and a unit square together with its center. We also prove that every two-coloring of three-dimensional space contains either a monochromatic unit square or a monochromatic noncoplanar tetrahedron with four edges of length one and two of length $\sqrt{2}$. Our proofs combine Székely's averaging inequalities for distance-avoiding sets with the Fourier-analytic and linear-programming methods of de Oliveira Filho and Vallentin, together with a theorem on rotation-invariant measures from OpenAI's recent work on the Hadwiger-Nelson problem.
\end{abstract}

\maketitle

\setcounter{tocdepth}{1}
\tableofcontents

\section{Introduction}

Euclidean Ramsey theory, initiated in a series of papers by Erd\H{o}s, Graham,
Montgomery, Rothschild, Spencer, and Straus, asks which finite configurations are unavoidable in
every coloring of Euclidean space \cite{Erdos1973,Erdos1975a,Erdos1975b}. We write $\mathbb E^d$ for $\mathbb R^d$ equipped with its Euclidean metric. For
two configurations $X,Y\subset\mathbb E^d$, we write $\mathbb E^d\to(X,Y)$ if
every red--blue coloring of $\mathbb E^d$ contains either a red $X$
or a blue $Y$ (up to isometry). Its negation, denoted
$\mathbb E^d\nrightarrow(X,Y)$, means that there exists a red--blue coloring with no
red $X$ and no blue $Y$. More generally, for configurations
$X_1,\ldots,X_r\subset\mathbb E^d$, write
$\mathbb E^d\to(X_1,\ldots,X_r)$ if every $r$-coloring of $\mathbb E^d$ contains a
configuration isometric to $X_i$ in some color $i$. Therefore
$\mathbb E^d\nrightarrow(X_1,\ldots,X_r)$ means that there exists an $r$-coloring
in which color $i$ contains no configuration isometric to $X_i$ for every $i$.

\subsection{Lines}

For $m \ge 2$, let
$\ell_m = \{0,1,\ldots,m-1\}\subset\mathbb R$ be the unit-spaced $m$-point
line. Erd\H{o}s, Graham, Montgomery, Rothschild, Spencer, and Straus proved
$\mathbb E^2\to(\ell_2,\ell_4)$ and also
$\mathbb E^2\to(\ell_2,K)$ for every three-point configuration $K\subset\mathbb E^2$
\cite{Erdos1975a}.

Subsequent work has extended these results in several directions. 
Tsaturian later advanced $\mathbb E^2\to(\ell_2,\ell_4)$ to
$\mathbb E^2\to(\ell_2,\ell_5)$ \cite{Tsaturian2017}. In dimension $3$,
Iv\'an established $\mathbb E^3\to(\ell_2,K)$ for every five-point configuration
$K\subset\mathbb E^3$, while Arman and Tsaturian showed
$\mathbb E^3\to(\ell_2,\ell_6)$ \cite{Ivan1979,Arman2018}.
An argument of Szlam gives, for some constant $c'>0$, $\mathbb E^d\to(\ell_2,K)$ for every
$K\subset\mathbb E^d$ with $|K| \le2^{c'd}$ \cite{Szlam2001,Conlon2019}. 

Our first result strengthens the lower bound (of Tsaturian) in the plane:

\begin{theorem}\label{thm:l2l7}
$\mathbb E^2\to(\ell_2,\ell_7)$.
\end{theorem}

In the opposite direction, Erd\H{o}s and Graham claimed, without proof, that
$\mathbb E^2\nrightarrow(\ell_2,\ell_m)$ for some $m<10^7$ \cite{ErdosGraham1979}. Conlon and Fox
proved
$\mathbb E^2\nrightarrow(\ell_2,\ell_{10^{10}})$ and, for
some constant $c>0$,
$\mathbb E^d\nrightarrow(\ell_2,\ell_m)$ whenever $d \ge 2$ and $m \ge 2^{cd}$
\cite{Conlon2019} by a nonconstructive argument. Recently, Currier, Mody, Xie, and
Zhang improved the planar bound to
$\mathbb E^2\nrightarrow(\ell_2,\ell_{6330})$. Their proof follows the
general outline of Conlon and Fox, using lattice-based colorings and
probabilistic arguments, and is also nonconstructive
\cite{CurrierModyXieZhang2026}. In forthcoming work, the second and third
authors of this paper show that $\mathbb E^2\nrightarrow(\ell_2,\ell_{42})$ by a computer-assisted explicit construction \cite{HuLiuToAppear}.
Combining these results leaves the cases $8 \le m \le41$ unresolved.

We next consider colorings with three or more colors. The Hadwiger--Nelson
problem asks for the chromatic number of the graph whose vertices are the
points of the plane and whose edges join pairs of points distance one apart. Such a coloring has no monochromatic $\ell_2$, so lower bounds for this chromatic
number imply corresponding results with repeated $\ell_2$. The Moser spindle shows
that this graph requires at least four colors, and hence
$\mathbb E^2\to(\ell_2,\ell_2,\ell_2)$ \cite{Moser1961}. As a corollary of de Grey's result that the chromatic number of the plane is at least five \cite{deGrey2018}, every four-coloring contains a monochromatic $\ell_2$; that is,
$\mathbb E^2\to(\ell_2,\ell_2,\ell_2,\ell_2)$.
OpenAI proved, using entirely new techniques, that the plane is not
five-colorable, raising the lower bound in the Hadwiger--Nelson problem to
six \cite{OpenAI}. Consequently,
$\mathbb E^2\to(\ell_2,\ell_2,\ell_2,\ell_2,\ell_2)$. For three colors,
Currier, Moore, and Yip proved that, for every $d \ge 2$,
$\mathbb E^d\nrightarrow(\ell_3,\ell_3,\ell_8)$,
$\mathbb E^d\nrightarrow(\ell_3,\ell_4,\ell_7)$, and
$\mathbb E^d\nrightarrow(\ell_3,\ell_5,\ell_5)$
\cite{CurrierMooreYipAvoiding}. Theorem~\ref{thm:l2l2l3} strengthens
the three-color conclusion from the Moser spindle by replacing one of the
configurations $\ell_2$ with $\ell_3$.

\begin{theorem}\label{thm:l2l2l3}
$\mathbb E^2\to(\ell_2,\ell_2,\ell_3)$.
\end{theorem}

We now return to two colors and consider lines in higher dimensions.
For two-color problems involving a red $\ell_3$, Erd\H{o}s, Graham, Montgomery, Rothschild,
Spencer, and Straus proved $\mathbb E^3\to(K,K)$ for every three-point configuration
$K\subset\mathbb E^3$, and hence $\mathbb E^3\to(\ell_3,\ell_3)$
\cite{Erdos1973}. In the plane, Currier, Moore, and Yip later proved $\mathbb E^2\to(\ell_3,\ell_3)$ with computer assistance
\cite{CurrierMooreYip3AP}. Conlon and
Wu proved $\mathbb E^d\nrightarrow(\ell_3,\ell_{10^{50}})$ for every $d \ge 2$
\cite{ConlonWu}. F\"uhrer and T\'oth improved $10^{50}$ to $1177$
\cite{FuhrerToth}, and Currier, Moore, and Yip later improved it to $20$
\cite{CurrierMooreYipAvoiding}. Thus, in the plane, there remains a gap $4 \le m \le19$. We were unable to improve the bounds in the plane, but in higher dimensions, we prove:

\begin{theorem}\label{thm:l3l4}
$\mathbb E^3\to(\ell_3,\ell_4)$.
\end{theorem}

\begin{theorem}\label{thm:l3l5}
$\mathbb E^4\to(\ell_3,\ell_5)$.
\end{theorem}

\begin{theorem}\label{thm:l3l6}
$\mathbb E^7\to(\ell_3,\ell_6)$.
\end{theorem}

Erd\H{o}s, Graham, Montgomery, Rothschild, Spencer, and Straus proved
$\mathbb E^d\nrightarrow(\ell_6,\ell_6)$ and
$\mathbb E^d\nrightarrow(\ell_4,\ell_4,\ell_4)$ for every $d \ge 2$
\cite{Erdos1973}. For two colors, interchanging the colors in
$\mathbb E^2\to(\ell_2,\ell_4)$ gives $\mathbb E^2\to(\ell_4,\ell_2)$
\cite{Erdos1975a}. Currier, Moore, and Yip proved
$\mathbb E^d\nrightarrow(\ell_4,\ell_{14})$ and
$\mathbb E^d\nrightarrow(\ell_5,\ell_8)$ for every $d \ge 2$
\cite{CurrierMooreYipAvoiding}. Thus the planar gap for red $\ell_4$ is
$3 \le m \le13$. In dimension $4$, we prove:

\begin{theorem}\label{thm:l4l4}
$\mathbb E^4\to(\ell_4,\ell_4)$.
\end{theorem}

Table~\ref{intro:tab:line-results} compares Theorems~\ref{thm:l2l7}--\ref{thm:l4l4}
with the results discussed above.

\begin{table}[htbp]
\centering
\small
\begin{tabular}{@{}llll@{}}
\toprule
Theorem & Prior lower bound & New lower bound & Upper bound\\
\midrule
\ref{thm:l2l7} & $\mathbb E^2\to(\ell_2,\ell_5)$ & $\mathbb E^2\to(\ell_2,\ell_7)$ & $\mathbb E^2\nrightarrow(\ell_2,\ell_{42})$\\
\ref{thm:l2l2l3} & $\mathbb E^2\to(\ell_2,\ell_2,\ell_2)$ & $\mathbb E^2\to(\ell_2,\ell_2,\ell_3)$ & \\
\addlinespace
\begin{tabular}[c]{@{}l@{}}
\ref{thm:l3l4}\\
\ref{thm:l3l5}\\
\ref{thm:l3l6}
\end{tabular}
& $\mathbb E^2\to(\ell_3,\ell_3)$
& \begin{tabular}[c]{@{}l@{}}
$\mathbb E^3\to(\ell_3,\ell_4)$\\
$\mathbb E^4\to(\ell_3,\ell_5)$\\
$\mathbb E^7\to(\ell_3,\ell_6)$
\end{tabular}
& $\mathbb E^d\nrightarrow(\ell_3,\ell_{20})$\\
\addlinespace
\ref{thm:l4l4} & $\mathbb E^2\to(\ell_4,\ell_2)$ & $\mathbb E^4\to(\ell_4,\ell_4)$ & $\mathbb E^d\nrightarrow(\ell_4,\ell_{14})$\\
\bottomrule
\end{tabular}
\caption{Line results.}
\label{intro:tab:line-results}
\end{table}

It is also worth mentioning a more general result of Conlon and F\"uhrer: if a finite configuration
$X\subset\mathbb E^d$ is not
contained in a sphere, then for some integer $m$,
$\mathbb E^n\nrightarrow(X,\ell_m)$ for every $n \ge d$
\cite{ConlonFuhrer2026}.

\subsection{Other configurations}

For $n \ge 3$ and $s>0$, write $W_n(s)$ for a regular $n$-gon of side length $s$ together with its center, and abbreviate $W_4(1)$ to $W_4$. Juh\'asz \cite{Juhasz1979} showed $\mathbb E^2\to(\ell_2,K)$ for every four-point subset $K\subset\mathbb E^2$, which immediately gives $\mathbb E^2\to(\ell_2,W_3(s))$ for all $s>0$, and constructed a $12$-point set $K\subset\mathbb E^2$ with $\mathbb E^2\nrightarrow(\ell_2,K)$. Csizmadia and T\'oth \cite{Csizmadia1994} cut twelve points down to eight with a centered heptagon of side $s = \frac95\sin(\pi/7)$, proving $\mathbb E^2\nrightarrow\bigl(\ell_2,W_7(\frac95\sin(\pi/7))\bigr)$. Among the open cases $4 \le n \le 6$, we settle $W_4$:

\begin{theorem}\label{thm:w4}
$\mathbb E^2\to(\ell_2,W_4)$.
\end{theorem}

Removing the center, let $P_n$ ($n \ge 3$) denote the regular $n$-gon of unit side length. $\mathbb E^2\to(\ell_2,P_4)$ is a corollary of Juh\'asz's theorem, and $\mathbb E^2\to(\ell_2,P_6)$ is immediate: if a red point $x$ exists, the six vertices of a unit regular hexagon centered at $x$ must all be blue, whereas if no point is red the entire plane is blue. For the regular pentagon $P_5$, we have:

\begin{theorem}\label{thm:p5}
$\mathbb E^2\to(\ell_2,P_5)$.
\end{theorem}

Let $e_1,e_2,e_3$ form an orthonormal basis of $\R^3$, and define
$\mathcal S = \{0,e_1,e_2,e_1+e_2\}$ and
$\mathcal T = \bigl\{0,e_1,e_2,\tfrac12e_1+\tfrac{\sqrt3}{2}e_3\bigr\}$.
These are a unit square and a noncoplanar tetrahedron. Each has four
pairwise distances equal to $1$ and two equal to $\sqrt2$.
A folklore construction gives $\mathbb E^2\nrightarrow(\mathcal S,\mathcal S)$
\cite{Gasarch2025}. Burr's argument gives
$\mathbb E^6\to(\mathcal S,\mathcal S)$ \cite{Erdos1973}, and
Cantwell lowered the dimension to $4$ \cite{Cantwell1996}.
We later show that every four-point set in
$\mathbb E^3$ with distance multiset $\{1,1,1,1,\sqrt2,\sqrt2\}$
is isometric to either $\mathcal S$ or $\mathcal T$. We prove:

\begin{theorem}\label{thm:square-tetrahedron}
Every two-coloring of $\mathbb E^3$ contains a monochromatic $\mathcal S$ or $\mathcal T$.
\end{theorem}

\subsection{Methods}

We briefly sketch a proof of $\mathbb E^2\to(\ell_2,\ell_7)$ (Theorem~\ref{thm:l2l7}) to illustrate our proof method. We first establish an inequality for colorings of a finite configuration, then average it and show that its expected value is negative. Let $P$ be the vertex set of an 18-point triangular lattice strip (see Figure~\ref{fig:l2l7-strip} in Section~\ref{sec:l2l7}). Let $r_v\in\{0,1\}$ indicate whether vertex $v\in P$ is red. We construct a quadratic polynomial $\mathcal{L}((r_v)_{v\in P})$ in these 18 variables that is nonnegative for every coloring of $P$ containing no red $\ell_2$ and no blue $\ell_7$. Its coefficients are chosen by linear programming. To extend this local inequality to colorings of the entire plane, we average it over all translates of the configuration $P$.

To average this inequality without assuming measurability, start with a coloring of $\mathbb{E}^2$ with no red $\ell_2$ and no blue $\ell_7$. Restrict it to $\F^2$, where $\F=\R\cap\overline{\mathbb Q}$ is the field of real algebraic numbers, and let $\mathcal{X}$ be the space of such colorings of $\F^2$. We show that $\mathcal{X}$ is compact and metrizable, and the group $G = \F^2\rtimes \mathrm{SO}(2,\F)$ acts on $\mathcal{X}$. Furthermore, this group is solvable and therefore amenable, so we can average our restricted coloring over a F\o lner sequence in $G$. Since $\mathcal{X}$ is compact, Prokhorov's theorem implies that a subsequence of these averages converges weakly to a probability measure $\mu$, which is invariant under $G$. Since (by construction) $\mathcal{L}$ depends only on a finite number of points, it is continuous on $\mathcal{X}$, and we can show that $\EE[\mathcal L] \ge 0$.

Now we aim to obtain an upper bound for this expected value. By translation
invariance, $\mathbb E_\mu[r_v] = \mathbb E_\mu[r_0]$ and
$\mathbb E_\mu[r_ur_v] = C(v-u)$, where $C(z) = \mathbb E_\mu[r_0r_z]$.
The correlation function is positive definite, so Bochner's theorem gives
\[
C(z) = \int_{\widehat{\mathbb F^2}}\chi(z)\,d\lambda(\chi),
\]
where $\lambda$ is a positive $\mathrm{SO}(2,\mathbb F)$-invariant measure. So instead of studying the original coloring, we can bound $\EE[\mathcal{L}]$ through its Fourier representation.

However, there is still some difficulty, as not every character of $\F^2$ extends continuously to $\R^2$. Recent work by \cite{OpenAI} shows that the part of $\lambda$ supported on these discontinuous characters is a multiple of Haar measure, whose Fourier coefficients vanish at every nonzero displacement. Therefore, we can simply focus on the continuous characters $z\mapsto e^{i\xi\cdot z}$, with $\xi\in\R^2$, whose rotational averages give Bessel integrals and allow us to bound the pair correlations appearing in $\EE[\mathcal L]$.

For our choice of coefficients, the resulting bound gives $\mathbb E_\mu[\mathcal L] \le-1$ using Sz\'ekely's averaging inequalities and the Fourier methods of de Oliveira Filho and Vallentin \cite{Szekely1984,deOliveira2010}. This contradicts $\mathbb E_\mu[\mathcal L] \ge 0$. Our main contribution is to construct local inequalities for the configurations considered here and apply them to arbitrary colorings (rather than measurable colorings). The other arguments use different finite configurations but follow a similar approach. In particular, the $W_4$ case requires an additional bound $\mathbb E_\mu[r_0] \le1/4$, and in dimensions $d \ge 3$ we use spherical averages.

\section{Fourier analysis of invariant measures}\label{sec:invariant-measures}

\subsection{Invariant probability measures on colorings}

Let $\F = \R\cap\overline{\mathbb Q}$. Since all finite configurations $P$ that appear in our local inequalities lie in $\F^d$, we can restrict any coloring of $\mathbb E^d$ avoiding the forbidden configurations to one on $\F^d$. For $k$ colors with forbidden configurations assigned to each $i\in\{1,\ldots,k\}$, let $\mathcal X\subseteq\{1,\ldots,k\}^{\F^d}$ be the set of valid colorings $\omega:\F^d\to\{1,\ldots,k\}$. In the product topology on $\{1,\ldots,k\}^{\F^d}$ (which is generated by cylinder sets fixing $\omega$ on finite subsets of $\F^d$), $\{1,\ldots,k\}^{\F^d}$ is compact by Tychonoff's theorem and metrizable via
\[
d(\omega,\eta) = \sum_{n = 1}^\infty 2^{-n}\mathbf 1_{\{\omega(p_n)\ne\eta(p_n)\}}\qquad(\F^d = \{p_1,p_2,\ldots\}).
\]
Note that the set $\mathcal X$ is closed. In fact, if $\omega\notin\mathcal X$, one of the forbidden configurations appears on a finite set $S\subset\F^d$. Every coloring agreeing with $\omega$ on $S$ has the same issue, so the complement of $\mathcal X$ is open, and $\mathcal X$ is both compact and metrizable.

Now let's use rotations. Define $m = \lfloor d/2\rfloor$ and $K = \{a+ib:a,b\in\F,\ a^2+b^2 = 1\}$, so an element $a+ib\in K$ rotates $\F^2$ by $(x,y)\mapsto(ax-by,bx+ay)$. The $j$th factor of $K^m$ acts on $\operatorname{span}(e_{2j+1},e_{2j+2})$ for $0 \le j<m$ and leaves the other coordinates fixed. For odd dimensions, it fixes the last coordinate as well. Together with translations, these rotations give $G_d = \F^d\rtimes K^m$, $(v,U)x = v+Ux$, and $(g\cdot\omega)(x) = \omega(g^{-1}x)$.

Since each $U\in K^m$ is orthogonal, $|(v+Ux)-(v+Uy)| = |x-y|$ for all $x,y\in\F^d$, so $G_d$ preserves $\mathcal X$. Since each $g\in G_d$ permutes the coordinates of $\{1,\ldots,k\}^{\F^d}$, the action is continuous. Now, we give $G_d$ the discrete topology. Since $G_d$ is an extension of the abelian group $\F^d$ by the abelian group $K^m$, it is amenable \cite[Theorem~G.2.1 and Proposition~G.2.2(ii)]{Bekka2008}. For a countable
discrete group, it is well known that amenability is equivalent to the F\o lner criterion: for
every finite $Q\subset G_d$ and every $\varepsilon>0$, there is a nonempty
finite $F\subset G_d$ with
$|gF\mathbin\triangle F|<\varepsilon|F|$ for all $g\in Q$
\cite[Theorem~G.5.1]{Bekka2008}. Enumerating $G_d = \{g_1,g_2,\ldots\}$ and applying this criterion to $Q = \{g_1,\ldots,g_n\}$ with $\varepsilon = 1/n$ gives a F\o lner sequence of finite sets $F_n\subset G_d$ satisfying $|gF_n\mathbin\triangle F_n|/|F_n|\to0$ as $n\to\infty$ for each $g\in G_d$.

\begin{lemma}\label{lem:law}
If $\mathcal X$ is nonempty, it supports a $G_d$-invariant Borel
probability measure $\mu$.
\end{lemma}
\begin{proof}
  Fix $\omega_0\in\mathcal X$ and let $\nu_n = |F_n|^{-1}\sum_{g\in F_n}\delta_{g\cdot\omega_0}$ be the law of $g\cdot\omega_0$ for $g$ drawn uniformly from $F_n$, where $\delta_\omega$ denotes the point mass at $\omega$. Thus $\nu_n(A)$ is the fraction of transformations in $F_n$ whose
resulting coloring belongs to the Borel set $A\subseteq\mathcal X$.

Probability measures on the compact metrizable space $\mathcal X$ are
weakly sequentially compact by Prokhorov's theorem
\cite[Theorem~5.1]{Billingsley1999}. After passing to a
subsequence, $\nu_n$
therefore converges weakly to a probability measure $\mu$. To prove
invariance, fix $g\in G_d$ and a continuous function $f$ on $\mathcal X$.
The two averages below have identical terms on $gF_n\cap F_n$, so
\begin{align*}
 \left|\int f(g\cdot\omega)\,d\nu_n(\omega)
          -\int f(\omega)\,d\nu_n(\omega)\right|& = 
\frac1{|F_n|}
 \left|\sum_{a\in gF_n}f(a\cdot\omega_0)
      -\sum_{a\in F_n}f(a\cdot\omega_0)\right|\\ 
      & \le \|f\|_\infty\frac{|gF_n\mathbin\triangle F_n|}{|F_n|}.
    \end{align*}
    As $n\to\infty$, the right-hand side vanishes by the F\o lner property. Since both $f$ and $\omega\mapsto f(g\cdot\omega)$ belong to $C(\mathcal X)$, weak convergence gives $\int f(g\cdot\omega)\,d\mu = \int f(\omega)\,d\mu$. Because $C(\mathcal X)$ separates Borel probability measures on $\mathcal X$, $\mu$ is $G_d$-invariant.
\end{proof}

Fix a $G_d$-invariant probability measure $\mu$ on $\mathcal X$ and write $\EE[f] = \int_{\mathcal X}f\,d\mu$. For $g\in G_d$, define $T_gf(\omega) = f(g^{-1}\cdot\omega)$ on $L^2(\mu)$, and write $T_z$ when $g = (z,I)$ is translation by $z\in\F^d$. Given $\phi:\{1,\ldots,k\}\to\R$, set $h(\omega) = \phi(\omega(0))$, so $T_zh(\omega) = \phi(\omega(z))$, and define the pair correlation $C_h(z) = \EE[hT_zh]$ for $z\in\F^d$ (when $h$ is the indicator of color $i$ at $0$, $C_h(z)$ is the probability that both $0$ and $z$ have color $i$). Since $\mu$ is translation-invariant, $\EE[(T_xh)(T_yh)] = C_h(y-x)$ for all $x,y\in\F^d$.

Recall that a function $F:\F^d\to\mathbb C$ is positive definite if $\sum_{i,j}\overline{c_i}c_jF(z_j-z_i) \ge 0$ for all $z_1,\ldots,z_n\in\F^d$ and $c_1,\ldots,c_n\in\mathbb C$. For $C_h$, we have
\begin{align*}
\sum_{i,j}\overline{c_i}c_jC_h(z_j-z_i)
& = \sum_{i,j}\overline{c_i}c_j\EE[(T_{z_i}h)(T_{z_j}h)]\\
& = \EE \left[ \left|\sum_i c_iT_{z_i}h\right|^2\right] \ge 0.
\end{align*}

Thus $C_h$ is positive definite, and we can apply Bochner's
theorem. Also, $C_h(0) = \EE[h^2]$, and invariance under
$K^m$ gives $C_h(Uz) = C_h(z)$ for every $U\in K^m$.

\subsection{Fourier representation in the plane}

For a countable discrete abelian group $A$, its Pontryagin dual $\widehat A = \operatorname{Hom}(A,\mathbb T)$ (where $\mathbb T = \{w\in\mathbb C:\vert{}w\vert{} = 1\}$) is a compact metrizable abelian group under the topology of pointwise convergence. Specialized to a countable discrete group $A$, Bochner's theorem for locally compact abelian groups \cite[Theorem~4.19]{Folland2016} takes the following form:

\begin{theorem}[Bochner]\label{thm:bochner}
Every positive-definite function $F:A\to\mathbb C$ on a countable discrete abelian group $A$ is the Fourier transform
\[
F(z) = \int_{\widehat A}\chi(z)\,d\lambda_F(\chi)\qquad(z\in A)
\]
of a unique positive Borel measure $\lambda_F$ on $\widehat A$, with $\lambda_F(\widehat A) = F(0)$.
\end{theorem}

Take $A = \F^2$ and set $D = \widehat{\F^2}$. Each $\xi\in\R^2$ defines a character $j(\xi)\in D$ via $j(\xi)(z) = e^{i\xi\cdot z}$, and $j:\R^2\to D$ is an injective map onto the subgroup $C = j(\R^2)$ of characters extending continuously to $\R^2$. Because $j$ maps closed balls in $\R^2$ homeomorphically onto compact subsets of $D$, the image $C$ is Borel and $j^{-1}:C\to\R^2$ is Borel measurable. Under the rotation action $(u\chi)(z) = \chi(uz)$ of $K$ on $D$, normalized Haar measure $m_D$ on $D$ satisfies $\int_D\chi(z)\,d m_D(\chi) = 0$ for all $z\in\F^2\setminus\{0\}$ (as evaluation at $z\ne0$ is a nontrivial character of $D$). By \cite[Theorem~2.3]{OpenAI}, $m_D$ is the unique $K$-invariant probability measure supported outside $C$:

\begin{theorem}\label{transfer:thm:openai}
Every $K$-invariant Borel probability measure on $D$ assigning measure zero to $C$ equals $m_D$.
\end{theorem}

For $\xi\in\R^2$, averaging $e^{i\xi\cdot z}$ over rotations gives the Bessel function $J_0(|\xi||z|)$, where
\[
J_0(t) = \frac1{2\pi}\int_0^{2\pi}e^{it\cos\theta}\,d\theta
\]
\cite[Eq.~10.9.1]{DLMF2026}. Note that if $z,w\in\F+i\F$ satisfy $|w| = |z|>0$, then $w/z\in K$, so $K$ acts transitively on each circle of radius $r>0$ in $\F^2$. Theorems~\ref{thm:bochner} and~\ref{transfer:thm:openai} give:

\begin{lemma}\label{lem:decomposition}
Every real-valued, $K$-invariant, positive-definite function $F:\F^2\to\R$ admits a representation
\begin{equation}\label{eq:decomposition}
F(z) = \int_{[0,\infty)}J_0(t|z|)\,d\nu(t)+a\one_{\{0\}}(z)\qquad(z\in\F^2),
\end{equation}
where $\nu$ is a finite nonnegative Borel measure on $[0,\infty)$, $a \ge 0$, and $\nu([0,\infty))+a = F(0)$.
\end{lemma}
\begin{proof}
Let $\lambda$ be the Bochner measure of $F$ on $D$; by uniqueness in Theorem~\ref{thm:bochner}, $\lambda$ and its restrictions $\lambda|_C$ and $\lambda|_{D\setminus C}$ are $K$-invariant. Put $a = \lambda(D\setminus C)$. If $a>0$, then $a^{-1}\lambda|_{D\setminus C} = m_D$ by Theorem~\ref{transfer:thm:openai}, so $\lambda|_{D\setminus C} = am_D$ for all $a \ge 0$, and integrating $\chi(z)$ against $\lambda|_{D\setminus C}$ gives $a\one_{\{0\}}(z)$.

On $C$, the pushforward $\rho = (j^{-1})_*(\lambda|_C)$ is a finite $K$-invariant Borel measure on $\R^2$. Since $K$ is dense in $\SO(2)$, dominated convergence implies that $\rho$ is $\SO(2)$-invariant. If $\sigma$ is normalized Haar measure on $\SO(2)$, Fubini's theorem gives
\[
\int_{\R^2}e^{i\xi\cdot z}\,d\rho(\xi)
 = \int_{\SO(2)}\int_{\R^2}e^{i(U\xi)\cdot z}\,d\rho(\xi)\,d\sigma(U)
 = \int_{\R^2}J_0(|\xi||z|)\,d\rho(\xi),
\]
and the pushforward $\nu$ of $\rho$ under $\xi\mapsto|\xi|$ satisfies \eqref{eq:decomposition} with $\nu([0,\infty)) = \lambda(C) = F(0)-a$.
\end{proof}

  Placing an atom of mass $a$ at $\infty\in[0,\infty]$ and setting $J_0(\infty\cdot|z|) = \one_{\{0\}}(z)$ folds $a\one_{\{0\}}(z)$ into \eqref{eq:decomposition}: for any real function $h$ of the color at $0$ under a $G_2$-invariant probability measure $\mu$, applying Lemma~\ref{lem:decomposition} to $C_h$ gives a nonnegative measure $\bar\nu_h$ on $[0,\infty]$ with
\begin{equation}\label{eq:raw}
C_h(z) = \int_{[0,\infty]}J_0(t|z|)\,d\bar\nu_h(t)\quad(z\in\F^2)
\end{equation}
and $\bar\nu_h([0,\infty]) = \EE[h^2]$.

\begin{corollary}\label{cor:centered}
Suppose $\mu$ is a $G_2$-invariant probability measure, and let $r$ be the indicator of a single color at the origin with mean $p = \EE[r]$. Then we can find a nonnegative Borel measure $\bar\eta$ on $[0,\infty]$ of total mass $\bar\eta([0,\infty]) = p(1-p)$ such that
\begin{equation}\label{eq:centered}
C_r(z) = p^2+\int_{[0,\infty]}J_0(t|z|)\,d\bar\eta(t)\qquad(z\in\F^2).
\end{equation}
\end{corollary}
\begin{proof}
We center the indicator by setting $h = r-p$. Because translation invariance gives $\EE[T_zr] = \EE[r] = p$, expanding the correlation directly gives $C_h(z) = \EE[(r-p)T_z(r-p)] = C_r(z)-p^2$, while at the origin we simply have $C_h(0) = \EE[h^2] = p(1-p)$. We then apply \eqref{eq:raw} to $h$ and set $\bar\eta = \bar\nu_h$ to obtain \eqref{eq:centered}.
\end{proof}

\subsection{Spherical kernels in higher dimensions}

Let $\sigma_{d-1}$ be normalized surface measure on $\mathbb S^{d-1}$ (where $d \ge 2$) and define $\Omega_d(t) = \int_{\mathbb S^{d-1}}e^{-it\omega_1}\,d\sigma_{d-1}(\omega)$ for $t\in\R$. By \cite[p.~1421]{deOliveira2010}, $\int_{\mathbb S^{d-1}}e^{ir\xi\cdot\omega}\,d\sigma_{d-1}(\omega) = \Omega_d(r|\xi|)$ for $\xi\in\R^d$ and $r \ge 0$. Clearly $\Omega_d$ is real and even with $\Omega_d(0) = 1$; on $[0,\infty]$ we put $\Omega_d(r\infty) = \one_{\{0\}}(r)$ ($r \ge 0$). In terms of Euler's gamma function $\Gamma$, the Bessel function $J_\alpha(t)$ ($\alpha,t \ge 0$) has the series expansion

\[
J_\alpha(t) = \sum_{n = 0}^{\infty}
\frac{(-1)^n}{n!\,\Gamma(n+\alpha+1)}
 \left(\frac t2\right)^{2n+\alpha}.
\]\begin{lemma}\label{hd:lem:omega}
For $d \ge 2$ and $t>0$,
\[
\Omega_d(t)
 = \Gamma \left(\frac d2\right)
 \left(\frac2t\right)^{d/2-1}J_{d/2-1}(t)
 = \sum_{n = 0}^{\infty}
\frac{(-1)^nt^{2n}}
{2^n n!\,d(d+2)\cdots(d+2n-2)},
\]
where the empty product is understood to be $1$. Moreover, for all $t\in\mathbb R$, we have the bound

\[
\vert{}\Omega_d'(t)\vert{} \le\frac{2\Gamma(d/2)}{(d-1)\sqrt\pi\,\Gamma((d-1)/2)}.
\]

\end{lemma}
\begin{proof}
We first recall the standard Bessel formula from \cite[Eq.~(1)]{deOliveira2010}. If we expand $J_{d/2-1}(t)$ as a power series and apply the identity $\Gamma(n+d/2) = 2^{-n}d(d+2)\cdots(d+2n-2)\Gamma(d/2)$, we can rewrite the $n$th term of $\Omega_d(t)$ as

\[
\Gamma \left(\frac{d}{2}\right) \left(\frac{2}{t}\right)^{d/2-1} \frac{(-1)^nt^{2n+d/2-1}}{2^{2n+d/2-1}n!\,\Gamma(n+d/2)} = \frac{(-1)^nt^{2n}}{2^n n!\prod_{j = 0}^{n-1}(d+2j)}.\]

To bound $\Omega_d'(t)$, note that for $\varepsilon\ne0$ the difference quotients satisfy
\[
 \left|\frac{e^{-i(t+\varepsilon)\omega_1}-e^{-it\omega_1}}{\varepsilon}\right| \le|\omega_1| \le1
\]
and converge pointwise to $-i\omega_1e^{-it\omega_1}$ as $\varepsilon\to0$; dominated convergence therefore gives $\Omega_d'(t) = -i\int_{\mathbb S^{d-1}}\omega_1e^{-it\omega_1}\,d\sigma_{d-1}(\omega)$ and hence $|\Omega_d'(t)| \le\int_{\mathbb S^{d-1}}|\omega_1|\,d\sigma_{d-1}(\omega)$.

Diaconis and Freedman \cite[Remark~2.10]{DiaconisFreedman1987}
give the density of the first $k$ coordinates of a uniformly distributed point on the sphere
of radius $\sqrt d$. Taking $k = 1$ and scaling the radius to one, the first
coordinate of a uniformly distributed point of $\mathbb S^{d-1}$ has
density
\[
\frac{\Gamma(d/2)}{\sqrt\pi\,\Gamma((d-1)/2)}
(1-u^2)^{(d-3)/2}
\]
on $[-1,1]$. Consequently,
\[
\begin{aligned}
\int_{\mathbb S^{d-1}}|\omega_1|\,\mathrm d\sigma_{d-1}(\omega)
& = \frac{2\Gamma(d/2)}{\sqrt\pi\,\Gamma((d-1)/2)}
\int_0^1u(1-u^2)^{(d-3)/2}\,\mathrm d u
\\
& = \frac{2\Gamma(d/2)}
{(d-1)\sqrt\pi\,\Gamma((d-1)/2)}.
\end{aligned}
\]
\end{proof}

For $d\in\{2,3,4,7\}$, Lemma~\ref{hd:lem:omega} together with the half-integer Bessel identities for $J_{1/2}$ and $J_{5/2}$ \cite[Eq.~(10.47.3) and Sec.~10.49(i)]{DLMF2026} specializes at $t>0$ (and by continuity at $t = 0$, where $\Omega_d(0) = 1$) to $\Omega_2(t) = J_0(t)$, $\Omega_3(t) = \frac{\sin t}{t}$, $\Omega_4(t) = \frac{2J_1(t)}{t}$, and $\Omega_7(t) = 15\frac{(3-t^2)\sin t-3t\cos t}{t^5}$.

\subsection{Fourier representation in higher dimensions}

To apply Theorem~\ref{transfer:thm:openai} one coordinate plane at a time, we need to introduce a lemma on the product of measures. 

\begin{lemma}\label{lem:haar-product}
Let $Y$ be a separable metric space and $\tau$ a finite positive Borel measure on $Y\times D$ that is invariant under $K$ acting on $D$. If $\tau(Y\times C) = 0$, then $\tau = \rho\otimes m_D$, where $\rho$ is the $Y$-marginal of $\tau$.
\end{lemma}
\begin{proof}
For each Borel set $B\subseteq Y$, the slice measure $\tau_B(A) = \tau(B\times A)$ on $D$ is $K$-invariant with $\tau_B(C) = 0$ and $\tau_B(D) = \rho(B)$. When $\rho(B)>0$, Theorem~\ref{transfer:thm:openai} gives $\rho(B)^{-1}\tau_B = m_D$, and when $\rho(B) = 0$ we have $\tau_B = 0$. Therefore we get $\tau(B\times A) = \rho(B)m_D(A)$ for all Borel rectangles $B\times A$, and thus $\tau = \rho\otimes m_D$.
\end{proof}

For $d \ge 2$, let $m = \lfloor d/2\rfloor$. Given $z\in\F^d$, let us write $z^{(i)} = (z_{2i-1},z_{2i})$ for $1 \le i \le m$, and set $\mathcal Z_d = \{z\in\F^d:z^{(i)} \ne 0\text{ for }1 \le i \le m\}.$

\begin{lemma}\label{lem:coordinate-planes}
Suppose $h$ is a real function of the color at $0$ under
a $G_d$-invariant probability measure $\mu$. Then there exists a positive Borel measure $\lambda_h$ of mass
at most $\EE[h^2]$, on $\R^d$ if $d = 2m$ and on
$\R^{2m}\times\widehat\F$ if $d = 2m+1$, such that
for every $z\in\mathcal Z_d$,
\begin{align}
C_h(z)
& = \int_{\R^d}e^{i\xi\cdot z}\,d\lambda_h(\xi)
&&\text{if }d = 2m,\label{eq:blocks-even}\\
C_h(v,s)
& = \int_{\R^{2m}\times\widehat\F}
e^{i\xi\cdot v}\chi(s)\,d\lambda_h(\xi,\chi)
&&\text{if }d = 2m+1,\label{eq:blocks-odd}
\end{align}
where $z = (v,s)$ in the odd-dimensional case.
\end{lemma}

\begin{proof}
Let $\lambda$ be the Bochner measure of $C_h$ on $D^m$ (if $d = 2m$) or $D^m\times\widehat\F$ (if $d = 2m+1$). Then $\lambda$ has mass $\EE[h^2]$ and is $K$-invariant in each of the $m$ factors of $D$. For $1 \le i \le m$, let $E_i$ be the set of tuples whose first $D$-coordinate outside $C$ is the $i$th one. The sets $E_1,\ldots,E_m$ are disjoint, cover the complement of $C^m$ (or $C^m\times\widehat\F$), and are $K$-invariant in the $i$th factor, so $\lambda|_{E_i} = \rho_i\otimes m_D$ by Lemma~\ref{lem:haar-product}. Since $z^{(i)}\ne0$, we have $\int_D\chi(z^{(i)})\,d m_D(\chi) = 0$, so $\int_{E_i}\chi(z)\,d\lambda = 0$. Let $\lambda_h$ denote the pushforward of $\lambda\vert{}_{C^m}$ (or $\lambda\vert{}_{C^m\times\widehat\F}$ when $d$ is odd) under $j^{-1}$ in each $C$ coordinate. It is immediate that $\lambda_h$ has total mass at most $\EE[h^2]$, and by construction it satisfies \eqref{eq:blocks-even}--\eqref{eq:blocks-odd}.
\end{proof}

When $d \ge 3$ is odd, it is convenient to write $z = (x,y,w,s)\in\F\times\F\times\F^{d-3}\times\F$ (where $w$ is simply empty if $d = 3$). Notice that whenever $y\ne0$, rotating the $(x,s)$-plane by any $U\in K$ leaves both $y$ and $w$ unchanged. Therefore, the $(x,y)$-projection of $Uz$ is still nonzero, and this allows us to apply \eqref{eq:blocks-odd} directly at $Uz$. With this in mind as motivation, we define $\mathcal Z_d^* = \{z\in\mathcal Z_d: z_2\ne0, (z_1,z_d)\ne(0,0)\}$, on which $K$ acts naturally by rotations in the coordinate plane spanned by $e_1$ and $e_d$.

\begin{lemma}\label{lem:odd-average}
Let $d \ge 3$ be odd, and let $h_1,\ldots,h_q$
be real functions of the color at $0$ under a
$G_d$-invariant probability measure $\mu$.
There are positive Borel measures $\Lambda_j$ on
$\R^d$, with $\Lambda_j(\R^d) \le\EE[h_j^2]$,
and a F\o lner sequence $\Phi_n\subset K$ such that
\begin{equation}\label{transfer:eq:nonzero-coordinate}
\lim_{n\to\infty}\frac1{|\Phi_n|}
\sum_{U\in\Phi_n}C_{h_j}(Uz)
 = \int_{\R^d}e^{i\xi\cdot z}\,d\Lambda_j(\xi)
\end{equation}
for every $z\in\mathcal Z_d^*$ and $1 \le j \le q$.
\end{lemma}

\begin{proof}
Let $\widetilde\lambda_j$ be the pushforward of $\lambda_{h_j}$
onto $\R^{d-2}\times D$ under $(\eta_x,\eta_y,\eta_w,\chi)\longmapsto
((\eta_y,\eta_w),\psi)$, $\psi(a,b) = e^{i\eta_xa}\chi(b).$ For $z = (x,y,w,s)$ in the lemma, define the bounded continuous function
\[
\Psi_z(\eta_y,\eta_w,\psi)
 = e^{i(\eta_y y+\eta_w\cdot w)}\psi(x,s).
\]
Set $A_U(\eta,\psi) = (\eta,\psi\circ U)$ for $U\in K$.
If $\Phi_n\subset K$ is a F\o lner sequence, we can put
\[
\tau_{j,n} = \frac1{|\Phi_n|}\sum_{U\in\Phi_n}
(A_U)_*\widetilde\lambda_j.
\]
Then use equation~\eqref{eq:blocks-odd} at $Uz$, to obtain
\[
\int\Psi_z\,d\tau_{j,n}
 = \frac1{|\Phi_n|}\sum_{U\in\Phi_n}C_{h_j}(Uz).
\]

Since $D$ is compact and the $\R^{d-2}$-marginal is independent of $n$, the measures $\tau_{j,n}$ are tight and have constant mass $M_j = \widetilde\lambda_j(\R^{d-2}\times D) \le\EE[h_j^2]$. Now we can use Prokhorov's theorem, which tells us that a common subsequence converges weakly to measures $\tau_j$, one for each $j$.

For $V\in K$ and $f\in C_b(\R^{d-2}\times D)$, we can exploit the F\o lner property to get
\[
 \left|\int f\circ A_V\,d\tau_{j,n}-\int f\,d\tau_{j,n}\right|
 \le M_j\|f\|_\infty\frac{|V\Phi_n\mathbin\triangle\Phi_n|}{|\Phi_n|}\longrightarrow0,
\]
so $\tau_j$ is $K$-invariant in the $D$ coordinate. Now, Lemma~\ref{lem:haar-product} shows that its restriction to $\R^{d-2}\times(D\setminus C)$ is $\rho_j\otimes m_D$, and importantly, this part has zero integral against $\Psi_z$ because $(x,s)\ne(0,0)$.

A character in $C$ has the form $\psi(a,b) = e^{i(\alpha a+\beta b)}$. Push the remaining part of $\tau_j$ forward by
$(\eta_y,\eta_w,\psi)\mapsto(\alpha,\eta_y,\eta_w,\beta)$, and call the resulting measure $\Lambda_j$. Note that its mass is at most $M_j$. The weak convergence above gives
\[
\lim_{n\to\infty}\int\Psi_z\,d\tau_{j,n}
 = \int\Psi_z\,d\tau_j
 = \int_{\R^d}e^{i\xi\cdot z}\,d\Lambda_j(\xi),
\]
as required.
\end{proof}

Let $\SO(d,\F)$ be the group of $d\times d$ orthogonal matrices of determinant $1$ with entries in $\F$. For $z_1,\ldots,z_N\in\F^d\setminus\{0\}$, real functions $h_1,\ldots,h_q$ of the color at $0$, and real coefficients $c,a_j,b_{j\ell}$, set $H_j(t) = a_j+\sum_{\ell = 1}^N b_{j\ell}\Omega_d(|z_\ell|t)$ ($t\in[0,\infty)$) and $H_j(\infty) = a_j$.

\begin{theorem}\label{transfer:thm:quadratic}
Let $d \ge 2$ and $\mu$ be a $G_d$-invariant probability measure. If
\begin{equation}\label{transfer:eq:hypothesis}
c+\sum_{j = 1}^q a_j\EE[h_j^2]
+\sum_{j = 1}^q\sum_{\ell = 1}^N b_{j\ell}C_{h_j}(Rz_\ell) \ge 0
\qquad\bigl(R\in\SO(d,\F)\bigr),
\end{equation}
then $c+\sum_{j = 1}^q\int_{[0,\infty]} H_j(t)\,d\nu_j(t) \ge 0$ for some nonnegative Borel measures $\nu_j$ on $[0,\infty]$ with $\nu_j([0,\infty]) = \EE[h_j^2]$.
\end{theorem}
\begin{proof}
For even $d$, let $\rho_j = \lambda_{h_j}$ from
Lemma~\ref{lem:coordinate-planes}. For odd $d$, we average
\eqref{transfer:eq:hypothesis} over $UR$, where $U\in\Phi_n$, and then take $\rho_j = \Lambda_j$ using Lemma~\ref{lem:odd-average}.
Either way, we have shown the existence of positive Borel measures $\rho_j$ on $\R^d$ of total mass $M_j: = \rho_j(\R^d) \le \EE[h_j^2]$ such that
\begin{equation}\label{transfer:eq:continuous-inequality}
Q(R) : = c+\sum_{j = 1}^q a_j\EE[h_j^2] +\sum_{j = 1}^q\sum_{\ell = 1}^N b_{j\ell} \int_{\R^d}\cos(\xi\cdot Rz_\ell)\,d\rho_j(\xi) \ge 0
\end{equation}
whenever $R\in\SO(d,\F)$ is chosen so that $Rz_1,\ldots,Rz_N$ satisfy the necessary conditions of the lemma.

Since $K$ is dense in $\SO(2)$ and coordinate-plane rotations generate $\SO(d)$, $\SO(d,\F)$ is dense in $\SO(d)$; the projection conditions fail only on a finite union of proper closed subvarieties of $\SO(d)$, so $Q \ge 0$ on a dense subset of $\SO(d)$. Continuity of $Q$ on $\SO(d)$ (by dominated convergence) then gives $Q(R) \ge 0$ on all of $\SO(d)$. Averaging $Q(R)$ over $R\in\SO(d)$ against normalized Haar measure $\sigma$ via Fubini's theorem gives
\[
\int_{\SO(d)}\int_{\R^d}\cos(\xi\cdot Rz_\ell)\,d\rho_j(\xi)\,d\sigma(R)
 = \int_{\R^d}\Omega_d(|\xi||z_\ell|)\,d\rho_j(\xi),
\]
and setting $\nu_j = \kappa_j+(\EE[h_j^2]-M_j)\delta_\infty$, where $\kappa_j$ is the pushforward of $\rho_j$ under $\xi\mapsto|\xi|$, gives $\nu_j([0,\infty]) = \EE[h_j^2]$ and $c+\sum_{j = 1}^q\int_{[0,\infty]} H_j\,d\nu_j = \int_{\SO(d)}Q\,d\sigma \ge 0$.
\end{proof}

\section{Two-colorings of the Plane}\label{sec:l2l7}
\begin{figure}[ht]
\centering
\begin{tikzpicture}[
  scale = .78,
  point/.style = {circle,fill = black,inner sep = 1.35pt},
  unit/.style = {black!55,line width = .45pt},
  column/.style = {black,line width = 1.2pt}
]
\foreach \i in {0,1}{
  \foreach \j in {0,1,2,3,4,5,6,7,8}{
    \coordinate (P\i\j) at ({\j+0.5*\i},{-0.8660254*\i});
  }
}
\draw[column] (P00)--(P01)--(P02)--(P03)--(P04)--(P05)--(P06)--(P07)--(P08);
\draw[column] (P10)--(P11)--(P12)--(P13)--(P14)--(P15)--(P16)--(P17)--(P18);
\foreach \j in {0,1,2,3,4,5,6,7,8}{
  \draw[unit] (P0\j)--(P1\j);
}
\draw[unit] (P01)--(P10);
\draw[unit] (P02)--(P11);
\draw[unit] (P03)--(P12);
\draw[unit] (P04)--(P13);
\draw[unit] (P05)--(P14);
\draw[unit] (P06)--(P15);
\draw[unit] (P07)--(P16);
\draw[unit] (P08)--(P17);
\foreach \i in {0,1}{
  \foreach \j in {0,1,2,3,4,5,6,7,8}{
    \node[point] at (P\i\j) {};
  }
}
\end{tikzpicture}
\caption{The eighteen-point triangular-lattice strip $P$. All drawn edges have length 1.}
\label{fig:l2l7-strip}

\end{figure}
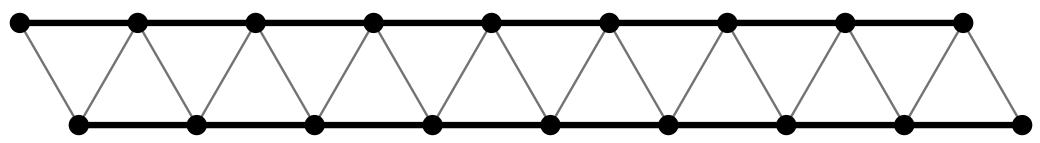

Let $P\subset\F^2$ be the eighteen-point triangular-lattice strip in Figure~\ref{fig:l2l7-strip}. Let $e_1 = (1,0)$ and $e_2 = (1/2,\sqrt3/2)$. For
$i\in\{0,1\}$ and $j\in\{0,1,\ldots,8\}$, put
$v_{ij} = (i-1/2)e_1+(j-4)e_2$, and set
$P = \{v_{ij}:i\in\{0,1\},\ j\in\{0,1,\ldots,8\}\}$. We use the label $ij$ for
$v_{ij}$. The formula gives $|v_{ij}| \le\sqrt{73}/2$ for every vertex.

For $i\in\{0,1\}$ and $s\in\{0,1,2\}$, put
$P_{i,s} = \{v_{ij}:s \le j \le s+6\}$. These are the six configurations
isometric to $\ell_7$ contained in $P$. The label difference
$(\Delta i,\Delta j)$ between consecutive points in a progression in the
triangular lattice satisfies
$\Delta i^2+\Delta i\Delta j+\Delta j^2 = 1$, whose six possibilities are
$(\pm1,0),(0,\pm1),(1,-1),(-1,1)$. Since $i\in\{0,1\}$, a progression of
seven points must have $\Delta i = 0$ and $|\Delta j| = 1$, so it is one of the
sets $P_{i,s}$.

For a coloring of $P$, let $r_p$ be the indicator that $p$ is red, let $c_p$ and $a_{pq} = a_{qp}$ be the integer coefficients listed in
Appendix~\ref{app:coefficient-data}, and define
\begin{equation}\label{l7:eq:local-functional}
\mathcal L = -1+\sum_{p\in P}c_pr_p
+\sum_{\{p,q\}\subset P}a_{pq}r_pr_q.
\end{equation}
When the coloring has no red $\ell_2$, every pair term at distance one
vanishes, and the remaining coefficients and indicators are integers, so the
local computation reduces to a simple integer calculation.

\begin{lemma}\label{l7:lem:local}
If a coloring of $P$ has no red $\ell_2$ and none of the six sets $P_{i,s}$
is entirely blue, then $\mathcal L \ge 0$.
\end{lemma}

\begin{proof}
Under the vertex ordering $00,01,\ldots,08,10,11,\ldots,18$, the $33$ pairs $\{p,q\}\subset P$ with $|v_p-v_q| = 1$ are
\[
\begin{aligned}
\{v_{ij},v_{i,j+1}\}&\quad(i\in\{0,1\},\ 0 \le j \le7),\\
\{v_{0j},v_{1j}\}&\quad(0 \le j \le8),\\
\{v_{0j},v_{1,j-1}\}&\quad(1 \le j \le8).
\end{aligned}
\]
Of the $2^{18}$ subsets $R\subseteq P$, $751$ of these contain none of the $33$ pairs and meet all six sets $P_{i,s}$. On each of the $751$ subsets, we can check that the integer expression \eqref{l7:eq:local-functional} is at least $0$ (Appendix~\ref{app:local-checks}).
\end{proof}

Suppose a coloring of $\mathbb E^2$ has no red $\ell_2$ and no blue $\ell_7$. Let $\mu$ be the $G_2$-invariant measure from Lemma~\ref{lem:law} and $r$ the red indicator at $0$; then $\EE[\mathcal L] \ge 0$ by Lemma~\ref{l7:lem:local}. Since $|e_1| = |e_2| = 1$ and $e_1\cdot e_2 = 1/2$, any two labels $p = (i,j)$ and $q = (i',j')$ satisfy $|v_p-v_q|^2 = Q(p-q)$ with $Q(a,b) = a^2+ab+b^2$; summing $a_{pq}$ over each value of $Q(p-q)$ (Table~\ref{l7:tab:distance-classes} in Appendix~\ref{l7:app:kernel-check}) defines the kernel
\begin{align}
K_7(t) = {}&-319-489J_0(t)+1035J_0(\sqrt3t)-669J_0(2t)\notag\\
&+600J_0(\sqrt7t)-700J_0(3t)+548J_0(\sqrt{13}t)\notag\\
&-355J_0(4t)+441J_0(\sqrt{21}t)-105J_0(5t)\notag\\
&+339J_0(\sqrt{31}t)-114J_0(6t)-139J_0(\sqrt{43}t)\notag\\
&+16J_0(7t)-90J_0(\sqrt{57}t).\label{l7:eq:HR}
\end{align}

\begin{proposition}\label{l7:prop:fourier}
There is a nonnegative Borel measure $\bar\nu_r$ on $[0,\infty]$ of mass $\EE[r]$ satisfying
\begin{equation}\label{l7:eq:fourier-error}
0 \le\EE[\mathcal L] = -1+\int_{[0,\infty]} K_7(t)\,d\bar\nu_r(t).
\end{equation}
\end{proposition}
\begin{proof}
We have $\sum_{p\in P}c_p = -319$ and $\EE[r_pr_q] = C_r(v_q-v_p)$. Applying \eqref{eq:raw} to $h = r$ (using $\EE[r^2] = \EE[r]$ and $J_0(0) = 1$) gives \eqref{l7:eq:fourier-error}.
\end{proof}

\begin{proposition}\label{l7:prop:kernel}
For all $t \ge 0$, $K_7(t)<0$.
\end{proposition}
\begin{proof}
See Appendix~\ref{l7:app:kernel-check}: adaptive stepping with $|K_7'| \le11682.22$ verifies $K_7(t)<0$ on $[0,100]$, while Lemma~\ref{shared:lem:analytic-bounds}(iv) gives $K_7(t)<-319+3000/\sqrt t \le-19$ for $t \ge 100$.
\end{proof}

\begin{proof}[Proof of Theorem~\ref{thm:l2l7}]
Because $K_7 \le0$ on $[0,\infty]$ (with $K_7(\infty) = -319$), Propositions~\ref{l7:prop:fourier} and~\ref{l7:prop:kernel} give $0 \le\EE[\mathcal L] \le-1$, which is absurd.
\end{proof}

\section{Three-colorings of the plane}\label{sec:l2l2l3}

\begin{figure}[ht]
\centering
\begin{tikzpicture}[
  scale = .82,
  point/.style = {circle,fill = black,inner sep = 1.35pt},
  edge/.style = {black!40,line width = .35pt}
]
\foreach \i in {0,1,2}{
  \foreach \j in {0,1,2}{
    \coordinate (P\i\j) at ({\i+0.5*\j},{0.8660254*\j});
  }
}
\draw[edge]
  (P00)--(P10)--(P20)
  (P01)--(P11)--(P21)
  (P02)--(P12)--(P22)
  (P00)--(P01)--(P02)
  (P10)--(P11)--(P12)
  (P20)--(P21)--(P22)
  (P10)--(P01)
  (P20)--(P11)--(P02)
  (P21)--(P12);
\foreach \i in {0,1,2}{
  \foreach \j in {0,1,2}{
    \node[point] at (P\i\j) {};
  }
}
\end{tikzpicture}
\caption{The nine-point triangular-lattice configuration $\Lambda$.}
\label{fig:l2l2l3-grid}
\end{figure}
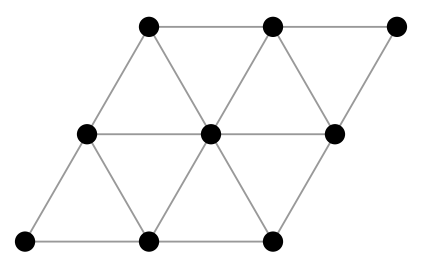

With $e_1 = (1,0)$ and $e_2 = (1/2,\sqrt3/2)$, put $v_{ij} = (i-1)e_1+(j-1)e_2$ ($i,j\in\{0,1,2\}$) and form the nine-point triangular grid $\Lambda = \{v_{ij}:i,j\in\{0,1,2\}\}$ of Figure~\ref{fig:l2l2l3-grid}, abbreviated by $ij$ for $v_{ij}$.

Using the same squared-distance form $Q(a,b) = a^2+ab+b^2$, we have $|v_{ij}-v_{i'j'}|^2 = Q((i,j)-(i',j'))$. The nonzero squared distances
in $\Lambda$ are $1,3,4,7,12$. The configurations isometric to $\ell_3$ used below lie in the
three unit directions $(1,0)$, $(0,1)$, and $(1,-1)$.

The four-element group $\mathcal G = \langle\sigma,\tau\rangle$ generated by $\sigma(i,j) = (j,i)$ and $\tau(i,j) = (2-i,2-j)$ acts isometrically on $\Lambda$, partitioning the $9$ vertices into $4$ orbits $\mathcal O_V$,
\[
\begin{array}{ccc}
\toprule
\text{representative}&\text{orbit}&\text{size}\\
\midrule
00&\{00,22\}&2\\
01&\{01,10,12,21\}&4\\
02&\{02,20\}&2\\
11&\{11\}&1\\
\bottomrule
\end{array}
\]
and the $\binom92 = 36$ unordered pairs into $13$ orbits $\mathcal O_P$. For $U\in\mathcal O_P$, write $q(U) = |v-w|^2$ ($\{v,w\}\in U$). Given $c:\Lambda\to\{G,R,B\}$, set $g = \one_G$ and $s = \one_R-\one_B$ (so $g^2 = g$, $s^2 = 1-g$, $gs = 0$), and define the orbit averages $[g]_O = |O|^{-1}\sum_{v\in O}g(v)$ ($O\in\mathcal O_V$) and
\[
[gg]_U = \frac1{|U|}\sum_{\{v,w\}\in U}g(v)g(w),\qquad
[ss]_U = \frac1{|U|}\sum_{\{v,w\}\in U}s(v)s(w)\qquad(U\in\mathcal O_P).
\]
With the integer coefficients $c_O,a_U,b_U$ listed in Tables~\ref{tab:c3-vertex-coefficients} and~\ref{tab:c3-pair-coefficients}, set
\[
\mathcal L(c) = -196+\sum_{O\in\mathcal O_V}c_O[g]_O+\sum_{U\in\mathcal O_P}\bigl(a_U[gg]_U+b_U[ss]_U\bigr).
\]
On the index set $V = \{0,1,2\}^2$, let
\[
\begin{aligned}
\mathcal E_1& = \{\{p,q\}\subseteq V:Q(p-q) = 1\},\\
\mathcal I_3& = \{(p,p+u,p+2u):u\in\{(1,0),(0,1),(1,-1)\},\ p,p+u,p+2u\in V\},
\end{aligned}
\]
and let $\mathcal A_3$ be the set of three-colorings $c:V\to\{G,R,B\}$ having no monochromatic red or blue pair in $\mathcal E_1$ and no all-green triple in $\mathcal I_3$.

\begin{lemma}\label{c3:lem:local}
Every coloring $c\in\mathcal A_3$ of $\Lambda$ with no red $\ell_2$, no blue $\ell_2$, and no green $\ell_3$ satisfies $\mathcal L(c) \ge 0$.
\end{lemma}

\begin{proof}
Since $|O|$ and $|U|$ divide $4$, the combination
\begin{align*}
4\mathcal L(c) = {}&-784+\sum_{O\in\mathcal O_V}\frac{4c_O}{|O|}\sum_{v\in O}g(v)\\
&+\sum_{U\in\mathcal O_P}\frac4{|U|} \left(a_U\sum_{\{v,w\}\in U}g(v)g(w)+b_U\sum_{\{v,w\}\in U}s(v)s(w)\right)
\end{align*}
is an integer for every $c$; checking all $|\mathcal A_3| = 304$ colorings in $\mathcal A_3$ (Appendix~\ref{app:local-checks}) gives $4\mathcal L(c) \ge 0$.
\end{proof}
Suppose a coloring of $\mathbb E^2$ has no red $\ell_2$, no blue $\ell_2$, and no green $\ell_3$, and let $\mu$ be the $G_2$-invariant measure from Lemma~\ref{lem:law} with $g = \one_G(0)$ and $s = \one_R(0)-\one_B(0)$, so $\EE[\mathcal L] \ge 0$. Collecting the coefficients in Tables~\ref{tab:c3-vertex-coefficients}--\ref{tab:c3-pair-coefficients} by $q(U)\in\{1,3,4,7,12\}$ gives the kernels
\begin{align}
H_g(t) = {}&-603-1050J_0(t)+1610J_0(\sqrt3t)-1152J_0(2t)\notag\\
&+1000J_0(\sqrt7t)+375J_0(\sqrt{12}t),\label{c3:eq:Hg}\\
H_s(t) = {}&-442J_0(t)+138J_0(\sqrt3t)-12J_0(2t)\notag\\
&+300J_0(\sqrt7t)+200J_0(\sqrt{12}t).\label{c3:eq:Hs}
\end{align}
\begin{proposition}\label{c3:prop:average}
There are nonnegative measures $\bar\nu_g,\bar\nu_s$ on $[0,\infty]$ with $\bar\nu_g([0,\infty])+\bar\nu_s([0,\infty]) = 1$ and
\begin{equation}\label{c3:eq:average}
0 \le\EE[\mathcal L] = -196+\int_{[0,\infty]} H_g\,d\bar\nu_g+\int_{[0,\infty]} H_s\,d\bar\nu_s.
\end{equation}
\end{proposition}
\begin{proof}
We have $\sum_{O\in\mathcal O_V}c_O = -603$ and $g^2 = g$, so applying \eqref{eq:raw} to $g$ and $s$ gives \eqref{c3:eq:average}, and $\bar\nu_g([0,\infty])+\bar\nu_s([0,\infty]) = \EE[g^2+s^2] = 1$.
\end{proof}
\begin{proposition}\label{c3:prop:kernels}
For all $t \ge 0$, $H_g(t)<193$ and $H_s(t)<193$.
\end{proposition}
\begin{proof}
Table~\ref{app:compact-bessel-checks} in Appendix~\ref{app:three-color-pentagon} bounds $H_g \le192.69$ and $H_s \le186.84$ on $[0,25]$ using grid step $1/400$, and Lemma~\ref{shared:lem:analytic-bounds}(iii) gives $H_g(t)<55$ and $H_s(t)<143$ for $t \ge 25$.
\end{proof}

\begin{proof}[Proof of Theorem~\ref{thm:l2l2l3}]
By Proposition~\ref{c3:prop:kernels}, as well as $H_g(\infty) = -603$ and $H_s(\infty) = 0$, we can show that $H_g,H_s \le193$ on $[0,\infty]$. Since $\bar\nu_g([0,\infty])+\bar\nu_s([0,\infty]) = 1$, we can apply Proposition~\ref{c3:prop:average} to get $0 \le\EE[\mathcal L] \le-196+193 = -3$, contradiction.
\end{proof}

\section{Lines in higher dimensions}\label{sec:hd-lines}

For a red--blue coloring, put $r = \one_R$ and
$s = \one_R-\one_B = 2r-1$. On a line $x+\R\omega$ with
$\omega\in\mathbb S^{d-1}$, write $r_j = r(x+j\omega)$ and
$s_j = s(x+j\omega)$. Avoiding a red $\ell_a$ and a blue $\ell_b$ simply means the binary sequence $(r_j)_{j\in\mathbb Z}$ can never contain $1^a$ or $0^b$ as a run of consecutive entries (and in the symmetric case $a = b = 4$, we never see four identical signs $s_j$ in a row). To handle Theorems~\ref{thm:l3l4}--\ref{thm:l4l4} all at once, we index the four problems by $q\in Q = \{1,2,3,4\}$ and record the ambient dimension $d_q$, the set of forbidden blocks $F_q$, and the step length $m_q$ in the table below:

\[
\begin{array}{ccccc}
\toprule
q & \text{relation} & d_q & F_q & m_q\\
\midrule
1 & (\ell_3,\ell_4) & 3 & \{111,0000\} & 3\\
2 & (\ell_3,\ell_5) & 4 & \{111,00000\} & 6\\
3 & (\ell_3,\ell_6) & 7 & \{111,000000\} & 6\\
4 & (\ell_4,\ell_4) & 4 & \{1111,0000\} & 4\\
\bottomrule
\end{array}
\]

Throughout what follows, we say a binary word avoids $F_q$ whenever none of the strings in $F_q$ appear inside it as a contiguous block. Define the de Bruijn vertex and edge sets (cf.\ Table~\ref{hd:tab:graph-sizes}) by $\mathcal V_q = \{u\in\{0,1\}^{m_q}:u\text{ avoids }F_q\}$ and
\[
\mathcal E_q = \bigl\{(u,v)\in\mathcal V_q^2:u_1\cdots u_{m_q-1} = v_0\cdots v_{m_q-2}\text{ and }u_0\cdots u_{m_q-1}v_{m_q-1}\text{ avoids }F_q\bigr\}.
\]
On the $(m_q+1)$-bit word $u_0\cdots u_{m_q-1}v_{m_q-1}$ of an edge $(u,v)\in\mathcal E_q$, put
\begin{align*}
w_{1}& = 4r_0-1-3r_0r_1-2r_0r_2,\\
w_{2}& = -6+28r_0-20r_0r_1-14r_0r_2-3r_0r_3-4r_0r_4+2r_0r_5+2r_0r_6,\\
w_{3}& = -2+11r_0-9r_0r_1-5r_0r_2-2r_0r_4-r_0r_5+r_0r_6,\\
w_{4}& = -1-6s_0s_1-4s_0s_2-2s_0s_3+s_0s_4,
\end{align*}
and take the integer potentials $\phi_q:\mathcal V_q\to\mathbb Z$ from Tables~\ref{hd:tab:phi34}--\ref{hd:tab:phi36} (Appendix~\ref{app:coefficient-data}).

\begin{proposition}\label{hd:prop:finite-inequalities}
For $q\in Q$ and $(u,v)\in\mathcal E_q$,
$w_q(u,v)+\phi_q(u)-\phi_q(v) \ge 0$.
\end{proposition}
\begin{proof}
We look at the coboundary difference $\Delta_q(u,v) = w_q(u,v)+\phi_q(u)-\phi_q(v)$ on each directed edge. Take, for instance, the edge $(000,001)\in\mathcal E_1$ corresponding to the four-letter word $0001$: here $w_1(0001) = -1$, and reading off $\phi_1(000) = 0$ and $\phi_1(001) = -1$ from Table~\ref{hd:tab:phi34} gives $\Delta_1(000,001) = -1+0-(-1) = 0$. Across $q = 1,2,3,4$, the graphs $\mathcal E_q$ have $12$, $73$, $78$, and $26$ edges, respectively. Once we plug in the integer potentials from Appendix~\ref{app:coefficient-data}, a quick check shows $\Delta_q(u,v) \ge 0$ across all edges (we record the full verification in Appendix~\ref{app:local-checks}).
\end{proof}

Next, take a $G_{d_q}$-invariant probability measure $\mu$ on colorings that avoid the monochromatic lines from row $q$. Whenever we fix a unit vector $\omega\in\F^{d_q}$, the colors at the points $0,\omega,\ldots,m_q\omega$ give us a valid directed edge $(u,v)\in\mathcal E_q$. Because $v$ is just a shift of $u$ by $\omega$, translation invariance gives $\EE[\phi_q(u)] = \EE[\phi_q(v)]$; telescoping the coboundary difference via Proposition~\ref{hd:prop:finite-inequalities} then forces $\EE[w_q] \ge 0$. For the asymmetric cases $q\in\{1,2,3\}$, we write $p = \EE[r]$ so that $\EE[rT_{k\omega}r] = p^2+C_{r-p}(k\omega)$ with variance $\EE[(r-p)^2] = p(1-p)$, whereas for $q = 4$ we simply use $s^2 = 1$. In Theorem~\ref{transfer:thm:quadratic}, the resulting spherical kernels are
\begin{align*}
H_{1}(t)& = -3\Omega_3(t)-2\Omega_3(2t),\\
H_{2}(t)& = -20\Omega_4(t)-14\Omega_4(2t)-3\Omega_4(3t)\\
&\phantom{ = {}}-4\Omega_4(4t)+2\Omega_4(5t)+2\Omega_4(6t),\\
H_{3}(t)& = -9\Omega_7(t)-5\Omega_7(2t)-2\Omega_7(4t)-\Omega_7(5t)+\Omega_7(6t),\\
H_{4}(t)& = -6\Omega_4(t)-4\Omega_4(2t)-2\Omega_4(3t)+\Omega_4(4t).
\end{align*}

\begin{proposition}\label{hd:prop:kernels}
For all $t \ge 0$, we have $H_{1}(t)<0.7$, $H_{2}(t)<2.95$, $H_{3}(t)<0.4$, and $H_{4}(t)<0.9$.
\end{proposition}
\begin{proof}
Table~\ref{hd:tab:kernel-checks} bounds each kernel on a finite interval. The tail bounds in Table~\ref{hd:tab:kernel-tails} cover the remaining values of $t$. Both verifications appear in Appendix~\ref{hd:app:line-kernels}.
\end{proof}

\begin{proof}[Proofs of Theorems~\ref{thm:l3l4}--\ref{thm:l4l4}]
Assume that a coloring avoids the monochromatic lines in the theorem
under consideration, and let $\mu$ be the corresponding invariant measure.
Theorem~\ref{transfer:thm:quadratic} gives the four inequalities
\begin{align*}
0& \le-1+4p-5p^2+\int H_{1}\,d\eta_{1},\\
0& \le-6+28p-37p^2+\int H_{2}\,d\eta_{2},\\
0& \le-2+11p-16p^2+\int H_{3}\,d\eta_{3},\\
0& \le-1+\int H_{4}\,d\nu_{4}.
\end{align*}

The first three measures have mass $p(1-p)$, while the fourth has mass one. Using the kernel bounds, which also hold at $\infty$ since $H_q(\infty) = 0$, we get the quadratic upper bounds as presented in Table~\ref{hd:tab:density-polynomials} for the first three cases. All three quadratics are negative for every real $p$, since their leading coefficients and discriminants are negative. In the fourth case, the upper bound is $-1+0.9 = -0.1$. Thus all four cases are impossible.

\begin{center}
\captionsetup{hypcap = false}
\begin{tabular}{@{}clc@{}}
\toprule
$q$ & Upper bound & Discriminant\\
\midrule
$1$ & $-1+4.7p-5.7p^2$ & $-0.71$\\
$2$ & $-6+30.95p-39.95p^2$ & $-0.8975$\\
$3$ & $-2+11.4p-16.4p^2$ & $-1.24$\\
\bottomrule
\end{tabular}
\captionof{table}{Polynomial bounds for the first three line inequalities.}\label{hd:tab:density-polynomials}
\end{center}
\end{proof}

\section{The centered unit square \texorpdfstring{$W_4$}{W4}}\label{sec:w4}

Figure~\ref{fig:planar-configurations} displays the twenty-point configuration of four copies of $W_4$ used to prove Theorem~\ref{thm:w4}, together with the nine-point configuration of two copies of $P_5$ used in Section~\ref{sec:p5}.

\begin{figure}[ht]
\centering
\begin{tikzpicture}[
  x = .75cm,
  y = .75cm,
  configvertex/.style = {circle,fill = black,inner sep = 1.15pt},
  configlabel/.style = {font = \scriptsize,inner sep = 1pt},
  pathedge/.style = {line width = .55pt},
  configsecond/.style = {line width = .55pt,dashed}
]
\begin{scope}[xshift = 0cm]
  \pgfmathsetmacro{\sdiag}{sqrt(2)/2}
  \pgfmathsetmacro{\slo}{1-\sdiag}
  \pgfmathsetmacro{\shi}{1+\sdiag}
  \draw[pathedge] ({-\sdiag},0)--(0,\sdiag)--(\sdiag,0)--(0,{-\sdiag})--cycle;
  \draw[pathedge] (\slo,0)--(1,\sdiag)--(\shi,0)--(1,{-\sdiag})--cycle;
  \draw[pathedge] (\slo,1)--(1,\shi)--(\shi,1)--(1,\slo)--cycle;
  \draw[pathedge] ({-\sdiag},1)--(0,\shi)--(\sdiag,1)--(0,\slo)--cycle;
  \foreach \p/\pos/\lab in {
    (0,0)/below left/$c_0$,(1,0)/below right/$c_1$,
    (1,1)/above right/$c_2$,(0,1)/above left/$c_3$}
    {\node[configvertex] at \p {};\node[configlabel,\pos] at \p {\lab};}
  \node[font = \small] at (.5,-1.9) {$(a)\;W_4$};
\end{scope}

\begin{scope}[xshift = 5cm]
  \coordinate (p50) at (0.85065,0);
  \coordinate (p51) at (0.26287,0.80902);
  \coordinate (p52) at (-0.68819,0.5);
  \coordinate (p53) at (-0.68819,-0.5);
  \coordinate (p54) at (0.26287,-0.80902);
  \coordinate (p55) at (-0.32492,1.61803);
  \coordinate (p56) at (-1.27598,1.30902);
  \coordinate (p57) at (-1.27598,0.30902);
  \coordinate (p58) at (-0.32492,0);
  \draw[pathedge] (p50)--(p51)--(p52)--(p53)--(p54)--cycle;
  \draw[configsecond] (p51)--(p55)--(p56)--(p57)--(p58)--cycle;
  \foreach \p/\pos/\lab in {
    p50/below right/0,p51/right/1,p52/above left/2,p53/below left/3,
    p54/below/4,p55/above right/5,p56/above left/6,p57/left/7,p58/below right/8}
    {\node[configvertex] at (\p) {};\node[configlabel,\pos] at (\p) {\lab};}
  \node[font = \small] at (-.2,-2.0) {$(b)\;P_5$};
\end{scope}

\end{tikzpicture}
\caption{The planar configurations used for $W_4$ and $P_5$.}
\label{fig:planar-configurations}
\end{figure}

\subsection{A bound on the red probability}
\label{sec:density}

For a $G_2$-invariant probability measure $\mu$ on $\mathcal X$, let $\mathcal H_c\subset L^2(\mu)$ be the subspace of functions $f$ such that $\|T_zf-f\|_{L^2(\mu)}\to0$ as $|z|\to0$ ($z\in\F^2$), and write $\lambda_f$ for the Bochner measure of $z\mapsto\EE[\overline f T_zf]$ on $D = \widehat{\F^2}$. By \cite[Lemma~4.1 and Proposition~4.2]{OpenAI}, $\mathcal H_c$ is the $L^2$-space of a $G_2$-invariant factor $\sigma$-algebra $\mathcal A$ and admits jointly measurable $\R^2$-translates:

\begin{theorem}\label{thm:continuous-functions}
There is a $G_2$-invariant sub-$\sigma$-algebra $\mathcal A$ such that $\mathcal H_c = L^2(\mathcal X,\mathcal A,\mu)$, with $f\in\mathcal H_c$ if and only if $\lambda_f(D\setminus C) = 0$, and $f\in\mathcal H_c^\perp$ if and only if $\lambda_f(C) = 0$; orthogonal projection onto $\mathcal H_c$ coincides with conditional expectation $\EE[\cdot\mid\mathcal A]$.
\end{theorem}

\begin{proposition}\label{prop:jointly-measurable}
  By $L^2$ continuity, the translation action of $\F^2$ on $\mathcal H_c$ extends, and for every $q\in\mathcal H_c$ with $0 \le q \le1$, there exists a jointly measurable function $H:\R^2\times\mathcal X\to[0,1]$ such that $H(x,\cdot) = T_xq$ in $L^2(\mu)$ for almost every $x\in\R^2$.
\end{proposition}

By Ambrus et al.\ \cite[Theorem~1]{Ambrus2024}, any measurable set $A\subset\R^2$ containing no unit-distance pair $\ell_2$ has upper density at most $0.2470<1/4$; writing $B_T = \{x\in\R^2:|x| \le T\}$, we use only the bound $1/4$:

\begin{theorem}\label{thm:measurable-density}
Every measurable set $A\subset\R^2$ with no $\ell_2$ satisfies $\limsup_{T\to\infty}|A\cap B_T|/|B_T| \le1/4$.
\end{theorem}

Following the density-point argument of Falconer \cite{Falconer1981} and Payne \cite{Payne2009} (cf.\ \cite[Lemma~4.3]{OpenAI}), call $x\in\R^2$ a density-one point of a measurable set $A\subset\R^2$ if $|A\cap B(x,\varepsilon)|/(\pi\varepsilon^2)\to1$ as $\varepsilon\to0$, where $B(x,\varepsilon)$ is the disk of radius $\varepsilon$ about $x$.

\begin{lemma}\label{lem:density-points}
Let $A\subset\R^2$ be measurable with $\one_A(x)\one_A(x+u) = 0$ for almost every $(x,u)\in\R^2\times\mathbb S^1$ under $\mathrm dx\,\mathrm d\sigma(u)$. Then the set $A^*$ of density-one points of $A$ is measurable, has $|A\mathbin\triangle A^*| = 0$, and contains no $\ell_2$.
\end{lemma}
\begin{proof}
The Lebesgue density theorem gives measurability of $A^*$ and $|A\mathbin\triangle A^*| = 0$. Suppose that $x_0,y_0\in A^*$ are distance one apart. For $0<\varepsilon<1/2$, consider $(x,x+u)$ such that $x\in B(x_0,\varepsilon)$, $u\in\mathbb S^1$, and $x+u\in B(y_0,\varepsilon)$. Finally, let $\nu_\varepsilon$ be the image of $\mathrm dx\,\mathrm d\sigma(u)$ on this set under $(x,u)\mapsto(x,x+u)$.

For $x\in B(x_0,\varepsilon/2)$, an arc of length at least $c_0\varepsilon$ of the circle $x+\mathbb S^1$ lies in $B(y_0,\varepsilon)$, so consequently
\[
\nu_\varepsilon\bigl(B(x_0,\varepsilon) \times B(y_0,\varepsilon)\bigr) \ge c\varepsilon^3
\]
for an absolute constant $c>0$. Additionally, each marginal of $\nu_\varepsilon$ has density at most $C\varepsilon$ with respect to planar Lebesgue measure, since for a fixed first or second coordinate, the admissible directions occupy an arc of length at most $C_0\varepsilon$.

Both $x_0$ and $y_0$ are density-one points of $A$. For sufficiently small $\varepsilon$,
\[
|B(x_0,\varepsilon)\setminus A|+|B(y_0,\varepsilon)\setminus A|
<\frac{c\varepsilon^2}{2C}.
\]
We can now use the marginal bounds to show that $\nu_\varepsilon(A\times A)>0$, which contradicts our assumption of almost-everywhere vanishing of $\one_A(x)\one_A(x+u)$.
\end{proof}

\begin{proposition}\label{prop:density}
Let $\mu$ be a $G_2$-invariant probability measure on colorings of $\F^2$ with no red $\ell_2$. If $r$ denotes the red indicator at $0$, then $\EE[r] \le1/4$.
\end{proposition}

\begin{proof}
Let $p = \EE[r]$ and $q = \EE[r\mid\mathcal A]\in\mathcal H_c$,
as in Theorem~\ref{thm:continuous-functions}, so that $0 \le q \le1$,
$\EE[q] = p$, and $r-q\in\mathcal H_c^\perp$. Both $q$ and
$r-q$ are $K$-invariant. Since the Bochner measure of $r-q$ has
$\lambda_{r-q}(C) = 0$, Theorem~\ref{transfer:thm:openai} gives
$\lambda_{r-q} = am_D$ for some $a \ge 0$; therefore
$\EE[(r-q)T_z(r-q)] = 0$ for $z\in\F^2\setminus\{0\}$.
Translations preserve $\mathcal H_c$ and its orthogonal complement,
and $\EE[rT_zr] = 0$ at unit distance. Hence
\[
0 = \EE[rT_zr]
 = \EE[qT_zq]+\EE[(r-q)T_z(r-q)]
 = \EE[qT_zq]\qquad(z\in\F^2,\ |z| = 1).
\]
It's now clear that by $L^2$-continuity, this identity holds with $z$ replaced by every
$u\in\mathbb S^1$.

Define $H:\R^2\times\mathcal X\to[0,1]$ as the jointly measurable
representative given by Proposition~\ref{prop:jointly-measurable}, and notice that $\EE[H(x,\cdot)] = p$ for almost every $x$, and by Fubini's theorem, $\EE[H(x,\cdot)H(x+u,\cdot)] = \EE[qT_uq] = 0$ for almost every
$(x,u)\in\R^2\times\mathbb S^1$. Now by Tonelli's theorem, for
$\mu$-almost every $\omega$ the set
$A_\omega = \{x\in\R^2:H(x,\omega)>0\}$ satisfies
Lemma~\ref{lem:density-points}, since its density-one representative
$A_\omega^*$ contains no $\ell_2$, and
$|A_\omega\mathbin\triangle A_\omega^*| = 0$.
As $H(\cdot,\omega) \le\one_{A_\omega}$, the bound applies to its ball averages.

For integers $n \ge 1$, set $Z_n(\omega) = |B_n|^{-1}\int_{B_n}H(x,\omega)\,d x$. Since $H(\cdot,\omega) \le\one_{A_\omega}$, Theorem~\ref{thm:measurable-density} implies $\limsup_n Z_n(\omega) \le1/4$ almost surely. On the other hand, we have $0 \le Z_n \le1$ and $\EE[Z_n] = p$. Now we apply Fatou's lemma to $1-Z_n$, which gives:
\[
1-p \ge \EE\bigl[\liminf_{n\to\infty}(1-Z_n)\bigr]
 = 1-\EE\bigl[\limsup_{n\to\infty}Z_n\bigr] \ge \frac{3}{4}
\]
and therefore $p \le 1/4$.
\end{proof}

\subsection{The local inequality for $W_4$}

Let $e_1 = (1,0)$, $e_2 = (0,1)$, and
\[
(u_0,u_1,u_2,u_3,u_4) = \frac{1}{\sqrt2}(0,e_1,e_2,-e_1,-e_2).
\]
For the unit-square vertices $c_0 = (0,0)$, $c_1 = (1,0)$, $c_2 = (1,1)$, and $c_3 = (0,1)$, set $v_{5k+j} = c_k+u_j$ for $0 \le k \le3$ and $0 \le j \le4$ (Figure~\ref{fig:planar-configurations}(a)). The four copies of $W_4$ have index sets $I_k = \{5k,\ldots,5k+4\}$. The $36$ unit-distance pairs are
\[
\mathcal E_4 = \{(i,j):0 \le i<j \le19,\ |v_i-v_j| = 1\},
\]
as listed in Appendix~\ref{app:w4-pairs}. If $r_i\in\{0,1\}$ is the red indicator at $v_i$, any $(\ell_2,W_4)$-free coloring of $\{v_0,\ldots,v_{19}\}$ satisfies $r_ir_j = 0$ for $(i,j)\in\mathcal E_4$, $\sum_{i\in I_k}r_i \ge 1$ for $0 \le k \le3$, and we set
\[
\mathcal L_{W_4}(r) = c+\sum_{i = 0}^{19}\alpha_i r_i+\sum_{0 \le i<j \le19}\beta_{ij}r_ir_j\qquad\bigl(r\in\{0,1\}^{20}\bigr),
\]
with $c = -24.994229$ and $10^6\alpha_i,10^6\beta_{ij}\in\mathbb Z$ from Tables~\ref{w4:tab:linear} and~\ref{w4:tab:coefficients}.

\begin{lemma}\label{w4:lem:local}
Every $r\in\{0,1\}^{20}$ with $r_ir_j = 0$ on $\mathcal E_4$ and $\sum_{i\in I_k}r_i \ge 1$ ($0 \le k \le3$) satisfies $\mathcal L_{W_4}(r) \ge 0$.
\end{lemma}
\begin{proof}
  We can check all $2338$ admissible vectors $r\in\{0,1\}^{20}$ to find that $10^6\mathcal L_{W_4}(r) \ge 0$ (see Appendix~\ref{app:local-checks}).
\end{proof}

Let $K_{W_4}(t) = \sum_{i<j}\beta_{ij}J_0(|v_i-v_j|t)$ (using Table~\ref{w4:tab:kernel}). Under a $G_2$-invariant measure $\mu$ on $(\ell_2,W_4)$-free colorings with $p = \EE[r]$, we have $C_r(z) = 0$ for $|z| = 1$, so Corollary~\ref{cor:centered} implies that $\int_{[0,\infty]}J_0(t)\,d\bar\eta(t) = -p^2$ and
\[
\int_{[0,\infty]} K_{W_4}(t)\,d\bar\eta(t)
 = \int_{[0,\infty]}\bigl(K_{W_4}(t)-600J_0(t)\bigr)\,d\bar\eta(t)-600p^2.
\]

\begin{proposition}\label{w4:prop:kernel}
For $t \ge 0$, $K_{W_4}(t)-600J_0(t)<46$.
\end{proposition}
\begin{proof}
The calculation in Appendix~\ref{w4:app:kernel-check} gives the bound on $[0,125]$, where we use Taylor estimates on $2500$ subintervals of length $1/20$. For $t \ge 125$, Lemma~\ref{shared:lem:analytic-bounds}(ii) gives a stronger estimate $|K_{W_4}(t)-600J_0(t)|<40$.
\end{proof}

\begin{proof}[Proof of Theorem~\ref{thm:w4}]
Let $\mu$ be the $G_2$-invariant measure obtained from an
$(\ell_2,W_4)$-free coloring. Proposition~\ref{prop:density}
gives $p = \EE[r] \le1/4$, and Corollary~\ref{cor:centered}
gives $\bar\eta([0,\infty]) = p(1-p)$.

The coefficient sums are $A = \sum_i\alpha_i = 139.366276$ and $B = \sum_{i<j}\beta_{ij} = 303.205864$. Use the local inequality, together with the estimate $K_{W_4}-600J_0 \le46$ (also at $\infty$), to obtain
\[
0 \le\EE[\mathcal L_{W_4}]
 \le c+(A+46)p+(B-646)p^2
 \le-24.99+185.37p-342.79p^2 = : \Phi(p).
\]
But $\Phi'(p) = 185.37-685.58p>13$ on $[0,1/4]$, and
\[
\boxed{\Phi(1/4) = -0.071875}.
\]
This forces $\Phi(p)<0$, which directly contradicts our first inequality.
\end{proof}

\section{The regular pentagon}\label{sec:p5}

To show that $\mathbb E^2\to(\ell_2,P_5)$, we consider a regular unit pentagon $B_0 = \{u_0, \ldots, u_4\}$ with circumradius $R_5 = 1/(2\sin(\pi/5))$ and vertices given explicitly by
\[u_j = \frac{1}{2\sin(\pi/5)} \left(\cos\frac{2\pi j}{5},\sin\frac{2\pi j}{5}\right)\]
for $0 \le j \le 4$.

Then, we translate $B_0$ by the edge vector $s = u_1-u_0$ to form a second pentagon $B_1 = B_0+s$. By construction, $B_0$ and $B_1$ meet at $u_1 = u_0+s$, and we claim they share no other vertices. To show this, note that any second overlap would require $u_i-u_j = s$. A unit chord of the regular pentagon is a side, and its ten directed side vectors are distinct; hence $u_i-u_j = u_1-u_0$ only for $(i,j) = (1,0)$, with indices modulo $5$.

Thus $B_0\cup B_1$ has nine points, which we will label by $v_i = u_i$ for $0 \le i \le4$ and $v_{4+j} = u_j+s$ for $1 \le j \le4$. It is now clear that two pentagons have index sets $I_0 = \{0,1,2,3,4\}$ and $I_1 = \{1,5,6,7,8\}$, and there are thirteen unit-distance pairs among these points. We will write $\mathcal E_5 = \{(i,j):i<j,\ |v_i-v_j| = 1\}$; Appendix~\ref{app:p5-pairs} lists the pairs and their distances.

Now, one can see that any $(\ell_2,P_5)$-free coloring of
$\{v_0,\ldots,v_8\}$ satisfies $r_ir_j = 0$ for $(i,j)  \in\mathcal E_5$, $\sum_{i\in I_0}r_i \ge 1$, and $\sum_{i\in I_1}r_i \ge 1$, where $r_i$ is the red indicator at $v_i$.
For $r\in\{0,1\}^9$, define
\[
\mathcal L_{P_5}(r) = -77+\sum_{i = 0}^8\alpha_i r_i
+\sum_{0 \le i<j \le8}\beta_{ij}r_ir_j.
\]
The coefficients are listed in Table~\ref{p5:tab:coefficients}
of Appendix~\ref{app:p5-pairs}.

\begin{lemma}\label{p5:lem:local}
If $r\in\{0,1\}^9$ satisfies $r_ir_j = 0$ for all $(i,j)\in\mathcal E_5$ and $\sum_{i\in I_0}r_i \ge 1$, $\sum_{i\in I_1}r_i \ge 1$, then $\mathcal L_{P_5}(r) \ge 0$.
\end{lemma}
\begin{proof}
Of the $2^9 = 512$ vectors $r\in\{0,1\}^9$, $42$ satisfy these constraints, and on each of those $42$ vectors, $\mathcal L_{P_5}(r)$ is at least $0$ (see Appendix~\ref{app:local-checks}).
\end{proof}

Collecting the pair coefficients $\beta_{ij}$ by $|v_i-v_j|^2$ gives
\begin{align*}
K_{P_5}(t) = {}&200J_0 \left(\sqrt{\frac{3-\sqrt5}{2}}\,t\right)-1048J_0(t)\\
&+100J_0 \left(\sqrt{\frac{5-\sqrt5}{2}}\,t\right)-385J_0 \left(\sqrt{\frac{3+\sqrt5}{2}}\,t\right).
\end{align*}

\begin{proposition}\label{p5:prop:kernel}
For all $t \ge 0$, $K_{P_5}(t)<393$.
\end{proposition}
\begin{proof}
In Appendix~\ref{app:three-color-pentagon} (Table~\ref{app:compact-bessel-checks}), evaluation at the $2400$ cell midpoints of $[0,12]$ with $|K_{P_5}'| \le1330$ gives $K_{P_5}(t) \le392.33$ on $[0,12]$, and Lemma~\ref{shared:lem:analytic-bounds}(ii) gives $|K_{P_5}(t)|<392$ for $t \ge 12$.
\end{proof}

\begin{proof}[Proof of Theorem~\ref{thm:p5}]
For an $(\ell_2,P_5)$-free coloring with $G_2$-invariant measure $\mu$ and $p = \EE[r]$, we have $\sum_i\alpha_i = 285$, $\sum_{i<j}\beta_{ij} = -1133$, and $K_{P_5} \le393$ on $[0,\infty]$ ($K_{P_5}(\infty) = 0$). Since $\bar\eta([0,\infty]) = p(1-p)$, Lemma~\ref{p5:lem:local} and Corollary~\ref{cor:centered} give
\[
\begin{aligned}
0 \le\EE[\mathcal L_{P_5}]& = -77+285p-1133p^2+\int_{[0,\infty]}K_{P_5}\,d\bar\eta\\
& \le-77+678p-1526p^2\\
& = -1-(1526p^2-678p+76) \le-1,
\end{aligned}
\]
which is impossible.
\end{proof}

\section{Unit squares and tetrahedra in \texorpdfstring{$\mathbb E^3$}{E3}}\label{sec:square-tetrahedron}

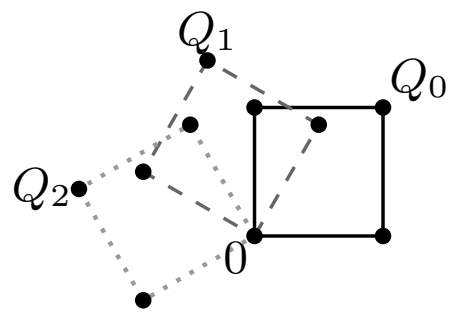
\begin{figure}[ht]
\centering
\begin{tikzpicture}[
  x = .85cm,
  y = .85cm,
  vertex/.style = {circle,fill = black,inner sep = 1.1pt},
  qzero/.style = {black,line width = .65pt},
  qone/.style = {black!60,dashed,line width = .65pt},
  qtwo/.style = {black!40,dotted,line width = .85pt},
  point label/.style = {font = \scriptsize,inner sep = 1pt}
]
\begin{scope}[xshift = 1.65cm]
  \coordinate (o) at (0,0);
  \coordinate (vzero) at (1,0);
  \coordinate (wzero) at (0,1);
  \coordinate (szero) at (1,1);
  \draw[qzero] (o)--(vzero)--(szero)--(wzero)--cycle;

  \coordinate (vone) at (.5,.866);
  \coordinate (wone) at (-.866,.5);
  \coordinate (sone) at (-.366,1.366);
  \draw[qone] (o)--(vone)--(sone)--(wone)--cycle;

  \coordinate (vtwo) at (-.5,.866);
  \coordinate (wtwo) at (-.866,-.5);
  \coordinate (stwo) at (-1.366,.366);
  \draw[qtwo] (o)--(vtwo)--(stwo)--(wtwo)--cycle;

  \foreach \p in {o,vzero,wzero,szero,vone,wone,sone,vtwo,wtwo,stwo}
    {\node[vertex] at (\p) {};}
  \node[point label,below left] at (o) {$0$};
  \node[point label,above right] at (szero) {$Q_0$};
  \node[point label,above] at (sone) {$Q_1$};
  \node[point label,left] at (stwo) {$Q_2$};
\end{scope}
\end{tikzpicture}
\caption{The three cross-sections $Q_r$ of the configuration
$\mathcal P$ in the plane perpendicular to the common edge $e$. Solid,
dashed, and dotted curves correspond to $Q_0$, $Q_1$, and $Q_2$, respectively.}
\label{fig:square-tetrahedron-configuration}
\end{figure}

Put $e = (1,0,0)\in\R^3$ and, for $r\in\{0,1,2\}$, let $v_r = \bigl(0,\cos(r\pi/3),\sin(r\pi/3)\bigr)$ and $w_r = \bigl(0,-\sin(r\pi/3),\cos(r\pi/3)\bigr).$ For each $r$, $(e,v_r,w_r)$ is an orthonormal basis of $\R^3$, so $Q_r = \{bv_r+cw_r:b,c\in\{0,1\}\}$ is a unit square in $e^\perp$ and $C_r = Q_r\cup(e+Q_r)$ is a unit cube. Since the six unit vectors $v_r,w_r$ and the three vectors $v_r+w_r$ of norm $\sqrt2$ have nine distinct directions in $e^\perp$, we have $Q_r\cap Q_s = \{0\}$ and $C_r\cap C_s = \{0,e\}$ ($r\ne s$), so $\mathcal P = C_0\cup C_1\cup C_2$ has $|\mathcal P| = 2+3\cdot6 = 20$ vertices (Figure~\ref{fig:square-tetrahedron-configuration}).

Both $\mathcal S$ and $\mathcal T$ have pairwise distance multiset
\[
\{1,1,1,1,\sqrt2,\sqrt2\}.
\]
The set $Y = \{0,w_0,e,w_1\}\subset\mathcal P$ has this multiset, and
\[
\det[w_0,e,w_1] = -\sqrt3/2\ne0.
\]
Thus $Y$ is noncoplanar and isometric to $\mathcal T$.

\begin{lemma}\label{sq:lem:four-point-classification}
Let $X\subset\mathbb E^3$ be a four-point set with distances
$\{1,1,1,1,\sqrt2,\sqrt2\}$. If the two $\sqrt2$-pairs are
disjoint, then $X\cong\mathcal S$, and otherwise
$X\cong\mathcal T$.
\end{lemma}

\begin{proof}
  First, suppose that $\{x_1,x_3\}$ and $\{x_2,x_4\}$ are disjoint $\sqrt2$-pairs, and translate the midpoint of $x_1x_3$ to $0$. Then $x_3 = -x_1$, $|x_1|^2 = 1/2$, and $|x_2\pm x_1| = |x_4\pm x_1| = 1$, so $x_1\perp x_2,x_4$ and $|x_2|^2 = |x_4|^2 = 1/2$. Since $|x_2-x_4|^2 = 2$, we get $x_4 = -x_2$ and $X = \{\pm x_1,\pm x_2\}\cong\mathcal S$.

If instead the $\sqrt2$-pairs share a vertex, we can translate the vertex belonging to neither pair to $0$ and write $X = \{0,a,b,c\}$ with $|a-b| = |b-c| = \sqrt2$ and $|a| = |b| = |c| = |a-c| = 1$. Then we can compute the Gram matrix of $(a,b,c)$ to get
\[
\begin{pmatrix}
1&0&1/2\\
0&1&0\\
1/2&0&1
\end{pmatrix},
\]
which is the same as the Gram matrix of the nonzero vertices of $\mathcal T$. Therefore $X\cong\mathcal T$.
\end{proof}

Consider the two isometries $\tau(x,y,z) = (1-x,y,z)$ and 
\[\sigma(x,y,z) = \left(x,-\frac{\sqrt3}{2}y-\frac{1}{2}z,-\frac{1}{2}y+\frac{\sqrt{3}}{2}z\right),\] 
and let $G = \langle\tau,\sigma\rangle$. We follow the labels $p_0,\ldots,p_{19}$ from Appendix~\ref{app:sq-labels}.

\begin{lemma}\label{sq:lem:configuration-counts}
The set $\mathcal P$ contains $22$ copies of $\mathcal S$ and $24$ copies of $\mathcal T$. The group $G$ preserves $\mathcal P$, has order $4$, and has $56$ orbits on its $190$ unordered pairs.
\end{lemma}
\begin{proof}
All points of $\mathcal P$ lie in $\mathbb Q(\sqrt3)^3$. Of the $\binom{20}{4} = 4845$ four-element subsets, the calculation in Appendix~\ref{sq:app:finite-check} finds $46$ with squared-distance multiset $\{1,1,1,1,2,2\}$. In $22$ of these the two pairs at squared distance $2$ are disjoint and in the other $24$ they meet, where we can use the preceding lemma to identify the two.

Now, we need to count the pair orbits. For $a,b,c\in\{0,1\}$, the map $\tau$ sends $ae+bv_r+cw_r$ to $(1-a)e+bv_r+cw_r$, so it preserves each cube $C_r$. Furthermore, $\sigma e = e$ and $(\sigma v_r,\sigma w_r) = (w_{2-r},v_{2-r})$, so $\sigma$ interchanges $C_0$ and $C_2$ and fixes $C_1$. The two maps are distinct commuting involutions, hence $|G| = 4$. Directly counting, their action on the $190$ pairs has $2$ orbits of size $1$, $14$ of size $2$, and $40$ of size $4$, giving us a total of $56$.
\end{proof}

With the $56$ pair orbits $\mathcal O_0,\ldots,\mathcal O_{55}$ and integer weights $a_0,\ldots,a_{55}$ ordered as in Table~\ref{sq:tab:pair-orbits}, set $X_j(s) = \sum_{\{p_u,p_v\}\in\mathcal O_j}s(p_u)s(p_v)$ for $s:\mathcal P\to\{-1,1\}$ and
\begin{equation}\label{sq:eq:local-functional}
\mathcal L_S(s) = 10^{-3} \left(-36608+\sum_{j = 0}^{55}a_jX_j(s)\right).
\end{equation}

\begin{lemma}\label{sq:lem:local}
Every $s:\mathcal P\to\{-1,1\}$ that is nonconstant on all $22$ copies of $\mathcal S$ and $24$ copies of $\mathcal T$ in $\mathcal P$ satisfies $\mathcal L_S(s) \ge 0$.
\end{lemma}
\begin{proof}
Testing all $2^{20}$ maps $s\in\{-1,1\}^{\mathcal P}$ leaves $47946$ valid colorings, and $10^3\mathcal L_S(s) \ge 0$ on each (Appendix~\ref{sq:app:finite-check}).
\end{proof}

For $\{p_u,p_v\}\in\mathcal O_j$, put $q_j = |p_u-p_v|^2$, sum $B_q = \sum_{0 \le j \le55:\,q_j = q}|\mathcal O_j|a_j$ (Table~\ref{sq:tab:distance-coefficients}; $\sum_q|B_q| = 620268$), and set
\begin{equation}\label{sq:eq:kernel}
K_S(t) = 10^{-3} \left(-36608+\sum_qB_q\Omega_3(\sqrt q\,t)\right).
\end{equation}

\begin{proposition}\label{sq:prop:average}
If a two-coloring of $\mathbb E^3$ avoids monochromatic copies of $\mathcal S$ and $\mathcal T$, then $\int_{[0,\infty]}K_S(t)\,d\nu(t) \ge 0$ for a Borel probability measure $\nu$ on $[0,\infty]$.
\end{proposition}
\begin{proof}
Apply Theorem~\ref{transfer:thm:quadratic} to $h = s = \one_R(0)-\one_B(0)$ ($s^2 = 1$) under the $G_3$-invariant measure $\mu$ of Lemma~\ref{lem:law}, using $\mathcal L_S \ge 0$ on $R\mathcal P$ ($R\in\SO(3,\F)$) from Lemma~\ref{sq:lem:local}.
\end{proof}

\begin{proposition}\label{sq:prop:kernel}
For all $t \ge 0$, $K_S(t)<-1/2$.
\end{proposition}
\begin{proof}
Appendix~\ref{sq:app:finite-check} checks $K_S(t)<-0.5225$ on $[0,40]$ via the grid $t_j = j/1000$ ($|K_S'|<931$), and bounds $K_S(t) \le-36.608+1240.536/t<-5.59$ for $t \ge 40$ via $|\Omega_3(\sqrt q\,t)| \le2/t$.
\end{proof}

\begin{proof}[Proof of Theorem~\ref{thm:square-tetrahedron}]
Because $K_S(t) \le-1/2$ for all $t\in[0,\infty]$ (with $K_S(\infty) = -36.608$), Proposition~\ref{sq:prop:average} gives $0 \le\int_{[0,\infty]}K_S\,d\nu \le-1/2$, which is impossible.
\end{proof}

\section{Linear programs for the local inequalities}\label{sec:lp}

For our planar two-color problems, we shall fix a finite point set $V\subset\F^2$ with unit-distance edges $E$ and target blue subsets $B_1,\ldots,B_s\subseteq V$. In this setting, a binary coloring $r = (r_v)_{v\in V}\in\{0,1\}^V$ is valid whenever $r_ur_v = 0$ for all $\{u,v\}\in E$ and $\sum_{v\in B_i}r_v \ge 1$ for $1 \le i \le s$. On such configurations, we study quadratic polynomials
\[
\widetilde{\mathcal L}(r) = c+\sum_{v\in V}a_vr_v+\sum_{\{u,v\}\subseteq V}b_{uv}r_ur_v
\]
with real coefficients $c,a_v$, and $b_{uv} = b_{vu}$.

Under a $G_2$-invariant measure $\mu$ with $p = \EE[r]$, Corollary~\ref{cor:centered} gives $C_r(u-v) = p^2+\int_{[0,\infty]}J_0(|u-v|t)\,d\bar\eta(t)$ with $\bar\eta([0,\infty]) = p(1-p)$. Now set $\widehat\eta = \bar\eta/(p(1-p))$ for $0<p<1$, and we have $\EE[\widetilde{\mathcal L}] = \int_{[0,\infty]}\widetilde H(p,t)\,d\widehat\eta(t)$, where
\[
\widetilde H(p,t) = c+p\sum_{v\in V}a_v+\sum_{\{u,v\}\subseteq V}b_{uv}\bigl(p^2+p(1-p)J_0(|u-v|t)\bigr).
\]
When $p = 0$ or $p = 1$, the function $\widetilde H(p,t)$ is independent of $t$, so any probability measure may be used in place of $\widehat\eta$.

For $W_4$, the unadjusted kernel $\widetilde H(p,t)$ is positive at $(p,t) = (1/4,0)$ ($\widetilde H(1/4,0) = 85.648806$), so the proof in Section~\ref{sec:w4} instead subtracts $600\int_{[0,\infty]}J_0(t)\,d\bar\eta(t) = -600p^2$ and combines the bound $K_{W_4}-600J_0<46$ with $p \le1/4$ (see Appendix~\ref{w4:app:kernel-check}).

For the eighteen-point strip $P$ of Section~\ref{sec:l2l7} (with valid coloring set $\mathcal A_7\subset\{0,1\}^P$) and the nine-point grid $\Lambda$ of Section~\ref{sec:l2l2l3}, the linear term is absorbed into the Bessel integral via \eqref{eq:raw}. On $P$, with variables $x_v$, $y_{uv} = y_{vu}$, and $\delta$, we put
\[
\begin{aligned}
\widetilde{\mathcal L}_7(r)& = -1+\sum_{v\in P}x_vr_v+\sum_{\{u,v\}\subseteq P}y_{uv}r_ur_v,\\
\widetilde K_7(t)& = \sum_{v\in P}x_v+\sum_{\{u,v\}\subseteq P}y_{uv}J_0(|u-v|t),
\end{aligned}
\]
and solve the LP (with $342$ bounds $|x_v|,|y_{uv}| \le100$ on the $171$ coefficients and grid $\mathcal N_7\subset[0,100]$; see Table~\ref{tab:lp-sizes}):
\[
\begin{array}{ll}
\operatorname{maximize}&\delta\\
\text{variables}&\delta,\ (x_v)_{v\in P},\ (y_{uv})_{\{u,v\}\subseteq P}\\
\operatorname{subject\ to}&\widetilde{\mathcal L}_7(r) \ge 0\quad(r\in\mathcal A_7),\\
&\widetilde K_7(t) \le-\delta\quad(t\in\mathcal N_7),\\
&|x_v|,|y_{uv}| \le100.
\end{array}
\]
On the nine-point grid $\Lambda$ with valid three-colorings $\mathcal A_3$ and orbit variables $(x_O)_{O\in\mathcal O_V}$, $(y_U)_{U\in\mathcal O_P}$, $(z_U)_{U\in\mathcal O_P}$, we set
\[
\widetilde{\mathcal L}_3 = -196+\sum_{O\in\mathcal O_V}x_O[g]_O+\sum_{U\in\mathcal O_P}\bigl(y_U[gg]_U+z_U[ss]_U\bigr)
\]
and
\[
\begin{aligned}
\widetilde H_g(t)& = \sum_{O\in\mathcal O_V}x_O+\sum_{U\in\mathcal O_P}y_UJ_0(\sqrt{q(U)}t),\\
\widetilde H_s(t)& = \sum_{U\in\mathcal O_P}z_UJ_0(\sqrt{q(U)}t),
\end{aligned}
\]
so $\EE[\widetilde{\mathcal L}_3] = -196+\int_{[0,\infty]}\widetilde H_g\,d\bar\nu_g+\int_{[0,\infty]}\widetilde H_s\,d\bar\nu_s$ with $\bar\nu_g([0,\infty])+\bar\nu_s([0,\infty]) = 1$. Imposing a common upper bound $196-\delta$ on $\widetilde H_g$ and $\widetilde H_s$ across $t\in\mathcal N_{\mathrm{col}}$ (with all $31$ variables unconstrained) gives the LP
\[
\begin{array}{ll}
\operatorname{maximize}&\delta\\
\text{variables}&\delta,\ (x_O)_{O\in\mathcal O_V},\ (y_U)_{U\in\mathcal O_P},\ (z_U)_{U\in\mathcal O_P}\\
\operatorname{subject\ to}&\widetilde{\mathcal L}_3(c) \ge 0\quad(c\in\mathcal A_3),\\
&\widetilde H_g(t) \le196-\delta,\ \widetilde H_s(t) \le196-\delta\quad(t\in\mathcal N_{\mathrm{col}}).
\end{array}
\]

In the four higher-dimensional line problems ($q\in\{1,2,3,4\}$) on the de Bruijn graphs $(\mathcal V_q,\mathcal E_q)$ of Section~\ref{sec:hd-lines}, write $\widetilde w_q(u,v)$ for the edge weight on $(u,v)\in\mathcal E_q$,
\[
\begin{aligned}
\widetilde w_{1}(r)& = -1+a_{1}r_0+c_{1,1}r_0r_1+c_{1,2}r_0r_2,\\
\widetilde w_{2}(r)& = -6+a_{2}r_0+\sum_{k = 1}^{6}c_{2,k}r_0r_k,\\
\widetilde w_{3}(r)& = -2+a_{3}r_0+\sum_{k\in\{1,2,4,5,6\}}c_{3,k}r_0r_k,\\
\widetilde w_{4}(s)& = -1+\sum_{k = 1}^{4}c_{4,k}s_0s_k,
\end{aligned}
\]
and $\widetilde H_q(t) = \sum_k c_{q,k}\Omega_{d_q}(kt)$ for the spherical kernel. For $q\in\{1,2,3\}$ (with $p = \EE[r]$ and $f = r-p$), Theorem~\ref{transfer:thm:quadratic} gives the (nonnegative) average
\[
\widetilde L_q = \alpha_q+\beta_qp + S_qp^2+\int_{[0,\infty]}\widetilde H_q(t) \,d\eta_f(t)
\]
where $\eta_f([0,\infty]) = p(1-p)$, $\beta_q = a_q$, $S_q = \sum_kc_{q,k}$, and $(\alpha_q,\beta_q,S_q)$ takes the rounded values $(-1,4,-5)$, $(-6,28,-37)$, and $(-2,11,-16)$ for $q = 1,2,3$; for $q = 4$, $\widetilde L_4 = -1+\int_{[0,\infty]}\widetilde H_4\,d\nu_4$ with $\nu_4([0,\infty]) = 1$. With target thresholds $(\kappa_1,\kappa_2,\kappa_3,\kappa_4) = (0.7,2.95,0.4,0.9)$ in dimensions $(d_1,d_2,d_3,d_4) = (3,4,7,4)$ and a reference vertex $u_q^*\in\mathcal V_q$ normalized by $\widetilde\phi_q(u_q^*) = 0$, the LP for case $q$ is
\[
\begin{array}{ll}
\operatorname{maximize}&\delta_q\\
\text{variables}&\delta_q,\ a_q,\ (c_{q,k})_k,\ (\widetilde\phi_q(u))_{u\in\mathcal V_q}\\
\operatorname{subject\ to}&\widetilde w_q(u,v)+\widetilde\phi_q(u)-\widetilde\phi_q(v) \ge 0\quad((u,v)\in\mathcal E_q),\\
&\widetilde H_q(t) \le\kappa_q-\delta_q\quad(t\in\mathcal N_q),\\
&\widetilde\phi_q(u_q^*) = 0,\ -100 \le\delta_q \le100,\ |a_q|,|c_{q,k}| \le100,\ |\widetilde\phi_q(u)| \le1000
\end{array}
\]

When $q = 4$, the variable $a_4$ is omitted. Table~\ref{tab:lp-sizes} lists the variable and constraint counts, and Appendix~\ref{app:checks} gives the sample grids $\mathcal N_q$. The rounded coefficients used in Section~\ref{sec:hd-lines} are
\[
\begin{array}{ccc}
\toprule
q&a_q&(c_{q,k})_k\\
\midrule
1&4&(-3,-2)\\
2&28&(-20,-14,-3,-4,2,2)\\
3&11&(-9,-5,-2,-1,1)\text{ for }k = 1,2,4,5,6\\
4&\text{--}&(-6,-4,-2,1)\\
\bottomrule
\end{array}
\]
and the integer functions $\phi_q$ are given in Tables~\ref{hd:tab:phi34}--\ref{hd:tab:phi36}.

For the $20$-point configuration $\mathcal P\subset\F^3$ of Section~\ref{sec:square-tetrahedron}, let $\mathcal A_{\mathcal S,\mathcal T}\subset\{-1,1\}^{\mathcal P}$ be the set of sign colorings avoiding monochromatic copies of the $22$ squares and $24$ tetrahedra of Lemma~\ref{sq:lem:configuration-counts}. With $X_j(s) = \sum_{\{p_u,p_v\}\in\mathcal O_j}s(p_u)s(p_v)$, set
\[
\begin{aligned}
\widetilde{\mathcal L}_{\mathcal S,\mathcal T}(s)& = c+\sum_{j = 0}^{55}b_jX_j(s),\\
\widetilde K_{\mathcal S,\mathcal T}(t)& = c+\sum_{j = 0}^{55}|\mathcal O_j|b_j\Omega_3(\sqrt{q_j}\,t).
\end{aligned}
\]
The LP (with solver bounds $|c|,|b_j| \le50$ and sample grid $\mathcal N_{\mathcal S,\mathcal T}$) is
\[
\begin{array}{ll}
\operatorname{maximize}&\delta\\
\text{variables}&\delta,c,(b_j)_{0 \le j \le55}\\
\operatorname{subject\ to}&\widetilde{\mathcal L}_{\mathcal S,\mathcal T}(s) \ge 0\quad\text{for }s\in\mathcal A_{\mathcal S,\mathcal T},\\
&\widetilde K_{\mathcal S,\mathcal T}(t) \le-\delta\quad\text{for }t\in\mathcal N_{\mathcal S,\mathcal T}.
\end{array}
\]
Its rounded solution $c = -36608/1000$ and $b_j = a_j/1000$ gives $\mathcal L_S(s)$ in \eqref{sq:eq:local-functional}. Table~\ref{tab:lp-sizes} summarizes the variable and constraint counts for all eight LPs.

\section{Finite subconfigurations and open problems}\label{sec:concluding}

In each of Theorems~\ref{thm:l2l7}--\ref{thm:square-tetrahedron}, the space $\mathcal X\subseteq\{1,\ldots,k\}^{\F^d}$ of valid colorings of $\F^d$ is empty, so compactness of $\{1,\ldots,k\}^{\F^d}$ reduces the theorem to a finite subset $V\subset\F^d$:

\begin{corollary}\label{cor:finite-radius}
For each of Theorems~\ref{thm:l2l7}--\ref{thm:square-tetrahedron}, the Ramsey property holds on a finite subset $V\subset\F^d$, and hence on every closed ball $B(0,R)\subset\mathbb E^d$ of sufficiently large radius $R$.
\end{corollary}

\begin{proof}
Let $\mathcal C_\alpha$ be the closed cylinder of colorings that
avoid a specified monochromatic occurrence on the finite vertex set
$V_\alpha\subset\F^d$. Then the space of admissible colorings is
$\mathcal X = \bigcap_{\alpha\in I}\mathcal C_\alpha = \emptyset$.
Because $\{1,\ldots,k\}^{\F^d}$ is compact, we can already find a finite subcollection $I_0\subset I$ with
\[
\bigcap_{\alpha\in I_0}\mathcal C_\alpha = \emptyset.
\]
If we now collect all the underlying vertices into the finite set $V = \bigcup_{\alpha\in I_0}V_\alpha$, then no $k$-coloring of $V$ can avoid our forbidden configurations, and it's clear that any Euclidean ball large enough to hold a translate of $V$ inherits the same property.
\end{proof}

Note that we can freely rescale $V$ by any $s>0$ to get the same result at scale $s$, or embed $V$ into $\mathbb E^n$ to handle any ambient dimension $n \ge d$.

\begin{corollary}\label{cor:infinitely-many}
  In each setting of Theorems~\ref{thm:l2l7}--\ref{thm:square-tetrahedron}, any coloring of $\mathbb E^d$ must contain infinitely many pairwise disjoint monochromatic copies of a single configuration.
\end{corollary}

\begin{proof}
  We can simply pack $\mathbb E^d$ with an infinite family of disjoint balls of radius $R$, where $R$ is chosen via Corollary~\ref{cor:finite-radius}. Each ball contains one of the finitely many forbidden monochromatic configurations. By the pigeonhole principle, a single color--configuration pair occurs in infinitely many of these balls, and the copies are disjoint.
\end{proof}

Because analogous LP searches did not give valid inequalities for $(\ell_2,\ell_8)$ in $\mathbb E^2$, for the line problems in lower dimensions, or for $(\ell_2,P_7)$ and $(\mathcal S,\mathcal S)$ in $\mathbb E^3$, we ask:

\begin{problem}
Does $\mathbb E^2\to(\ell_2,\ell_8)$ hold?

\end{problem}

\begin{problem}
Can any of the dimensions in Theorems~\ref{thm:l3l4}--\ref{thm:l4l4}
be lowered?
\end{problem}

\begin{problem}
Does $\mathbb E^3\to(\mathcal S,\mathcal S)$ hold?
\end{problem}

\begin{problem}
Does $\mathbb E^2\to(\ell_2,P_7)$ hold?
\end{problem}

\begin{problem}
For each relation in Corollary~\ref{cor:finite-radius}, what is the
infimum of the radii of balls on which it holds?
\end{problem}

\section*{Acknowledgments}

This work was conducted as part of the REU on Combinatorics, AI, and Algorithms for Real Problems (REU-CAAR) at the University of Maryland, with support from the National Science Foundation under Grant No. 2447466. We are grateful to William Gasarch for introducing us to this problem and for his guidance throughout the project. We also thank Charlie Klawitter, Eva Wu, Hridhaan Banerjee, and Tane Park for helpful discussions, as well as Kelin Zhu, Neal Shandilya, and Sandy Yin for proofreading an earlier draft. Finally, we acknowledge the use of OpenAI's GPT-5.6 Luna for assistance with proofreading, implementing the linear programs, and other computational aspects of this work.

\appendix
\setlength{\tabcolsep}{3pt}
\renewcommand{\arraystretch}{0.88}
\captionsetup[table]{skip = 3pt}

\section{Computational verifications}\label{app:checks}

The code used for the computational verifications is available at \url{https://github.com/pretzeledkoala/measurable-lines}.

Numerical estimates use intervals rounded outwards. Exhaustive coloring checks
use the listed integers or exact terminating decimals.

\subsection{Linear-program sampling data}

The kernel sample sets used in the linear programs are
\[
\begin{aligned}
\mathcal N_7& = \{j/400:0 \le j \le40000\},&
\mathcal N_{\mathrm{col}}& = \{j/400:0 \le j \le10000\},\\
\mathcal N_{1}& = \{j/5000:0 \le j \le30000\},&
\mathcal N_{2}& = \{j/2000:0 \le j \le50000\},\\
\mathcal N_{3}& = \{j/10000:0 \le j \le100000\},&
\mathcal N_{4}& = \{j/2500:0 \le j \le62500\},\\
\mathcal N_{\mathcal S,\mathcal T}& = \{j/1000:0 \le j \le40000\}
\end{aligned}
\]
where $j\in\mathbb Z$. For $(\ell_2,P_5)$, the LP uses $2140$ pairs $(p,t)$ to select the coefficients, and Appendix~\ref{app:three-color-pentagon} checks the resulting kernel $K_{P_5}$; for $W_4$, Appendix~\ref{w4:app:kernel-check} checks $[0,125]$ on $2500$ intervals of width $1/20$; and for $\mathcal S$--$\mathcal T$, the LP includes box bounds $|c|,|b_j| \le50$.

Using $J_0(x) = \pi^{-1}\int_0^\pi\cos(x\sin\theta)\,d\theta$ \cite[Eq.~10.9.1]{DLMF2026} and Hankel's expansion with remainder \cite[Eq.~10.17.3 and Sec.~10.17(iii)]{DLMF2026}, we record the analytic bounds used below.

\begin{lemma}\label{shared:lem:analytic-bounds}\,
\begin{enumerate}
  \item[(i)] $|J_0'(x)| < 2/3$ for $x\in\mathbb R$.
  \item[(ii)] $|J_0(x)| \le \sqrt{2/(\pi x)}$ for $x \ge 2$.
  \item[(iii)] $|J_0(\sqrt{q} t)| < 1/(6 q^{1/4})$ for $q \ge 1$ and $t \ge 25$.
  \item[(iv)] $|J_0(\sqrt{q} t)| < 5/(6\sqrt t\,q^{1/4})$ for $q \ge 1$ and $t \ge 100$.
  \item[(v)] $|\Omega_3(t)| \le 1/t$ for $t>0$.
  \item[(vi)] $|\Omega_4(t)| \le 2/t$ for $t>0$.
  \item[(vii)] $|\Omega_7(t)| \le 30/t^3$ for $t \ge 10$.
\end{enumerate}
\end{lemma}
\begin{proof}
From $J_0(x) = \pi^{-1}\int_0^\pi\cos(x\sin\theta)\,d\theta$, differentiation under the integral sign gives
\[
|J_0'(x)| \le\frac{1}{\pi}\int_{0}^{\pi}\sin{\theta}\,d\theta = \frac2\pi<\frac{2}{3}.
\]
For (ii), set $\chi = x-\pi/4$ and write
\[
J_0(x) = \sqrt{\frac2{\pi x}}\bigl(A(x)\cos\chi+B(x)\sin\chi\bigr),
\]
where $A = 1-9/(128x^2)+R_e$ and $B = 1/(8x)-R_o$ with $0 \le R_e \le0.113x^{-4}$ and $0 \le R_o \le0.074x^{-3}$. When $x \ge 2$, we have $0<1-9/(128x^2) \le A \le1-0.04x^{-2}$ (as $0.113/4<9/128-0.04$) and $0 \le B \le1/(8x)$ (as $0.074<0.125x^2$), so
\[
A^2+B^2 \le(1-0.04x^{-2})^2+0.015625x^{-2} = 1-0.064375x^{-2}+0.0016x^{-4} \le1,
\]
and Cauchy--Schwarz gives (ii). Putting $x = \sqrt q\,t$ in (ii) with $\pi>3$ gives (iii)--(iv); $\Omega_3(t) = (\sin t)/t$ gives (v); $|J_1(t)| \le1$ \cite[Eq.~10.9.2]{DLMF2026} gives (vi); and $|\Omega_7(t)| \le15(t^2+3t+3)/t^5 \le30/t^3$ ($t \ge 10$) gives (vii).
\end{proof}

For $z\in\R$ and $u = z^2/4$, truncating $J_0(z) = \sum_{m = 0}^\infty(-1)^mu^m/(m!)^2$ at $m = N$ when $u<(N+2)^2$ gives
\begin{equation}\label{shared:eq:bessel-tail}
 \left|\sum_{m = N+1}^{\infty}\frac{(-1)^m u^m}{(m!)^2}\right|
 \le\frac{u^{N+1}}{((N+1)!)^2}\frac1{1-u/(N+2)^2}.
\end{equation}
When $q,t\in\mathbb Q_{ \ge 0}$ and $z = \sqrt q\,t$, one can show that the
parameter $u = qt^2/4$, the $N$th partial sum, and \eqref{shared:eq:bessel-tail}
are all rational. Let $I$ be an interval on which $|K'| \le D$. Suppose that
for every $t\in I$ there is a sample point $s\in I$ with $|t-s| \le h/2$ and
$K(s) \le M$. Then the mean value theorem gives $K(t) \le M+Dh/2$ for every
$t\in I$.

\subsection{Local coloring enumerations}\label{app:local-checks}

On the $18$-point strip $P$, setting $r_p(R) = \one_{p\in R}$ gives $F(R) = -1+\sum_{p\in P}c_pr_p(R)+\sum_{\{p,q\}\subset P}a_{pq}r_p(R)r_q(R) \ge 0$ on all $751$ subsets $R\subseteq P$ avoiding the $33$ unit edges and meeting all six lines $P_{i,s}$. On the $9$-point grid $\Lambda$ ($g = \one_G$, $s = \one_R-\one_B$), we have $4\mathcal L \ge 0$ on all $304$ valid three-colorings; for $W_4$ and $P_5$, we have $10^6\mathcal L_{W_4} \ge 0$ on all $2338$ valid colorings and $\mathcal L_{P_5} \ge 0$ on all $42$ valid colorings. For the four line graphs $(\mathcal V_q,\mathcal E_q)$ ($m_q = 3,6,6,4$):
\begin{table}[H]
\centering
\small
\begin{tabular}{cccc}
\toprule
\text{case}&\text{excluded consecutive patterns}&\text{vertices}&\text{edges}\\
\midrule
$(\ell_3,\ell_4)$&111,\ 0000&7&12\\
$(\ell_3,\ell_5)$&111,\ 00000&41&73\\
$(\ell_3,\ell_6)$&111,\ 000000&43&78\\
$(\ell_4,\ell_4)$&1111,\ 0000&14&26\\
\bottomrule
\end{tabular}
\caption{Vertex and edge counts for the four line graphs.}
\label{hd:tab:graph-sizes}
\end{table}
\noindent and the potentials $\phi_q$ of Appendix~\ref{app:coefficient-data} satisfy $w_q(u,v)+\phi_q(u)-\phi_q(v) \ge 0$ on all $|\mathcal E_q|$ edges.

\subsection{Verification of the $W_4$ kernel bound}\label{w4:app:kernel-check}

Define $H(t) = K_{W_4}(t)-600J_0(t) = \sum_q c_qJ_0(\sqrt q\,t)$ using the coefficients from Table~\ref{w4:tab:kernel} (after
subtracting $600$ from $c_1$), and then partition $[0,125]$ into
$2500$ intervals $I$ of width $h = 1/20$ and midpoint $m_I$.
We can now compute
$u_I \ge H(m_I)$ and $d_I \ge |H'(m_I)|$, and truncate each
$J_0$ or $J_1$ series only when the bound
\eqref{shared:eq:bessel-tail} is below $10^{-80}$.
Since $|J_0''| \le1$, we can apply Taylor's theorem to get
\[
H(t) \le u_I+\frac h2d_I+\frac{h^2}{8}\sum_q|c_q|q<46
\qquad(t\in I).
\]
This estimate holds on all $2500$ intervals.

For the tail, we have $\sum_q|c_q|q^{-1/4}<542$. By Lemma~\ref{shared:lem:analytic-bounds}(ii) and $\pi>3$,
\[
|H(t)| \le\sqrt{\frac{2}{3t}}\sum_q|c_q|q^{-1/4}<40
\qquad(t \ge 125).
\]
We now use the same coefficient tables to get $\sum_i\alpha_i = 139.366276$ and $\sum_{i<j}\beta_{ij} = 303.205864$, hence $A+46 = 185.366276$ and $B-646 = -342.794136$.

\subsection{Verification of the $K_7$ kernel bound}\label{l7:app:kernel-check}

Table~\ref{l7:tab:distance-classes} lists the pair count $m_q$ and coefficient sum $A_q$ for each squared distance in
\[
\mathcal Q = \{1,3,4,7,9,13,16,21,25,31,36,43,49,57\},
\]
where $K_7(t) = -319+\sum_{q\in\mathcal Q}A_qJ_0(\sqrt q\,t)$; the classes $q = 64,73$ have $A_q = 0$.
\begin{table}[ht]
\centering
\small
\begin{tabular}{crr}
\toprule
$q$&$m_q$&$A_q$\\
\midrule
$1$&$33$&$-489$\\
$3$&$15$&$1035$\\
$4$&$14$&$-669$\\
$7$&$13$&$600$\\
$9$&$12$&$-700$\\
$13$&$11$&$548$\\
$16$&$10$&$-355$\\
$21$&$9$&$441$\\
$25$&$8$&$-105$\\
$31$&$7$&$339$\\
$36$&$6$&$-114$\\
$43$&$5$&$-139$\\
$49$&$4$&$16$\\
$57$&$3$&$-90$\\
\bottomrule
\end{tabular}
\caption{Distance classes and coefficient sums for $K_7$.}
\label{l7:tab:distance-classes}
\end{table}
Using $|K_7'| \le D: = 11682.216002$ (from Lemma~\ref{shared:lem:analytic-bounds}(i)), we step across $0 = t_0<t_1<\cdots<t_L = 100$: given an upper bound $K_i^+ \ge K_7(t_i)$ with $K_i^+<0$, choosing $t_{i+1}$ with $D(t_{i+1}-t_i) \le-K_i^+/2$ ensures $K_7(t) \le K_i^++D(t-t_i) \le K_i^+/2<0$ on $[t_i,t_{i+1}]$, which reaches $t_L = 100$. For $t \ge 100$, $\sum_{q\in\mathcal Q}|A_q|q^{-1/4}<3600$ and Lemma~\ref{shared:lem:analytic-bounds}(iv) give $K_7(t)<-319+3000/\sqrt t \le-19$.

\subsection{Verification of the three-color and pentagon kernels}\label{app:three-color-pentagon}

With grid spacing $h$, sample maximum $M$, and derivative bound $|K'| \le D$ (from Lemma~\ref{shared:lem:analytic-bounds}(i)), the bound $M+Dh/2$ on bounded intervals is summarized in Table~\ref{app:compact-bessel-checks}:
\begin{table}[H]
\centering
\footnotesize
\begin{tabular}{@{}ccrrrr@{}}
\toprule
kernel&interval&$h$&$M$&$D$&interval bound\\
\midrule
$H_g$&$[0,25]$&$0.0025$&$180$&$10151$&$192.69$\\
$H_s$&$[0,25]$&$0.0025$&$184.076$&$2208$&$186.84$\\
$K_{P_5}$&$[0,12]$&$0.005$&$389$&$1330$&$392.33$\\
\bottomrule
\end{tabular}
\caption{Sample, derivative, and interpolation bounds for $H_g$, $H_s$, and $K_{P_5}$.}
\label{app:compact-bessel-checks}
\end{table}
For $H_g$ and $H_s$ (sampled at $t_j = j/400\in[0,25]$), Table~\ref{app:compact-bessel-checks} gives $H_g,H_s<193$ on $[0,25]$; for $t \ge 25$, Lemma~\ref{shared:lem:analytic-bounds}(iii) with $3^{-1/4}<0.77$, $4^{-1/4}<0.715$, $7^{-1/4}<0.625$, $12^{-1/4}<0.556$ gives $H_g(t)<55$ and $H_s(t)<143$. For $K_{P_5}$ (sampled at the $2400$ cell midpoints $t_j = (2j+1)/400\in[0,12]$), Table~\ref{app:compact-bessel-checks} gives $K_{P_5}(t)<393$ on $[0,12]$; for $t \ge 12$, writing $(q_1,q_2,q_3,q_4) = ((3-\sqrt5)/2,1,(5-\sqrt5)/2,(3+\sqrt5)/2)$ and $(c_1,c_2,c_3,c_4) = (200,-1048,100,-385)$, the bounds $q_1^{1/4}>0.78$, $q_3^{1/4}>1.08$, $q_4^{1/4}>1.27$ give $\sum_{i = 1}^4|c_i|q_i^{-1/4}<1701$, whence $|K_{P_5}(t)|<8505/\sqrt{471}<392$ by Lemma~\ref{shared:lem:analytic-bounds}(ii) (as $\sqrt{q_1}t>2$ and $\pi>3.14$).

\subsection{Verification of the four line kernels}\label{hd:app:line-kernels}

From Lemma~\ref{hd:lem:omega}, we have $|\Omega_3'| \le1/2$, $|\Omega_4'| \le4/(3\pi)$, and $|\Omega_7'| \le5/16$, which give the derivative bounds $D$ and mean-value bounds $M+Dh/2$ in Table~\ref{hd:tab:kernel-checks} (where each grid of step $h$ includes both endpoints):
\begin{table}[ht]
\centering
\footnotesize
\begin{tabular}{@{}ccrrrr@{}}
\toprule
kernel&interval&$h$&$M$&$D$&interval bound\\
\midrule
$H_{1}$&$[0,6]$&$0.0002$&$0.69$&$3.5$&$0.6904$\\
$H_{2}$&$[0,25]$&$0.0005$&$2.938$&$47.5$&$2.9499$\\
$H_{3}$&$[0,10]$&$0.0001$&$0.399$&$11.875$&$0.3996$\\
$H_{4}$&$[0,25]$&$0.0004$&$0.887$&$12$&$0.8894$\\
\bottomrule
\end{tabular}
\caption{Sample, derivative, and interpolation bounds for the four line kernels.}\label{hd:tab:kernel-checks}
\end{table}
\noindent On the complementary rays, Lemma~\ref{shared:lem:analytic-bounds}(v)--(vii) bounds $|c_k\Omega_3(kt)| \le|c_k|/(kt)$, $|c_k\Omega_4(kt)| \le2|c_k|/(kt)$, and $|c_k\Omega_7(kt)| \le30|c_k|/(kt)^3$ ($t \ge 10$), giving Table~\ref{hd:tab:kernel-tails} and completing the proof of Proposition~\ref{hd:prop:kernels}.
\begin{table}[H]
\centering
\footnotesize
\begin{tabular}{@{}cccc@{}}
\toprule
kernel&range&tail bound&bound on the range\\
\midrule
$H_{1}$&$t \ge 6$&$4/t$&$0.667$\\
$H_{2}$&$t \ge 25$&$60/t$&$2.4$\\
$H_{3}$&$t \ge 10$&$291/t^3$&$0.291$\\
$H_{4}$&$t \ge 25$&$18/t$&$0.72$\\
\bottomrule
\end{tabular}
\caption{Tail bounds for the four line kernels.}
\label{hd:tab:kernel-tails}
\end{table}

\subsection{Computational verification for the square--tetrahedron configuration}
\label{sq:app:finite-check}

The geometric and local counts in Section~\ref{sec:square-tetrahedron} follow
from arithmetic in $\mathbb Q(\sqrt3)$. Represent $a+b\sqrt3$ by $(a,b)$, with
addition and multiplication given by $(a,b)+(c,d) = (a+c,b+d)$ and $(a,b)(c,d) = (ac+3bd,ad+bc).$
We now compute $d_{ij} = |p_i-p_j|^2$ between the $20$ points $p_0,\ldots,p_{19}$ (by brute force). Among the $\binom{20}{4} = 4845$ four-point subsets, $46$ have squared distances $\{1,1,1,1,2,2\}$. Of these, $22$ are copies of $\mathcal S$ and $24$ are copies of $\mathcal T$, depending on whether the two longest edges are disjoint or share a vertex. 

The group $G = \langle\tau,\sigma\rangle$ has four elements. Its action on the
unordered pairs is obtained from
\[
\mathcal O(i,j) = 
\bigl\{\{\min(g(i),g(j)),\max(g(i),g(j))\}:g\in G\bigr\},
\]
and a calculation gives the following orbit-size distribution:
\[
\begin{array}{ccccc}
\toprule
\text{orbit size}&1&2&4&\text{total}\\
\midrule
\text{number of orbits}&2&14&40&56\\
\text{pairs accounted for}&2&28&160&190
\\
\bottomrule
\end{array}
\]
The representatives and the associated distances and coefficients are listed
in Appendix~\ref{app:coefficient-data}.

Substituting the integer coefficients from Appendix~\ref{app:coefficient-data}
into the local expression gives $1000\mathcal L_S \ge 0$ for every allowed
coloring.

On the grid $t_j = j/1000$ in $[0,40]$, interval evaluation gives
$K_S(t_j)<-0.988$. Every squared distance in the kernel is less than $9$,
and Lemma~\ref{hd:lem:omega} gives $|\Omega_3'| \le1/2$, so
$|K_S'(t)|<931$. Each point of $[0,40]$ is within $1/2000$ of a grid
point; hence $K_S(t)<-0.988+0.4655 = -0.5225<-1/2$.

For $t \ge 40$, the smallest squared distance is $2-\sqrt3>1/4$, so
$|\Omega_3(\sqrt q\,t)| \le2/t$. Therefore
$K_S(t) \le-36.608+1240.536/t \le-5.5946<-1/2$.

\subsection{Linear-program sizes}

Table~\ref{tab:lp-sizes} records the dimensions and scalar constraint counts
of the eight linear programs used in the coefficient searches. The $W_4$ local
and kernel verifications are recorded separately in this appendix. The columns
labeled local, kernel, and variable bounds count scalar
inequalities. The equality column counts scalar equalities. The total number
of constraints is the sum of these four columns. A two-sided variable bound
gives two scalar inequalities.

\begin{table}[ht]
\centering
\scriptsize
\begin{tabular}{@{}lrrrrrr@{}}
\toprule
LP&variables&local&kernel&equalities&variable bounds&total constraints\\
\midrule
$(\ell_2,\ell_7)$
  &$172 = 1+18+153$&751&40001&0&342&41094\\
$(\ell_2,\ell_2,\ell_3)$
  &$31 = 1+4+13+13$&304&20002&0&0&20306\\
$(\ell_3,\ell_4)$
  &$11 = 1+1+2+7$&12&30001&1&22&30036\\
$(\ell_3,\ell_5)$
  &$49 = 1+1+6+41$&73&50001&1&98&50173\\
$(\ell_3,\ell_6)$
  &$50 = 1+1+5+43$&78&100001&1&100&100180\\
$(\ell_4,\ell_4)$
  &$19 = 1+4+14$&26&62501&1&38&62566\\
$(\ell_2,P_5)$
  &$47 = 2+9+36$&42&2140&0&94&2276\\
$\mathcal S$--$\mathcal T$
  &$58 = 1+1+56$&47946&40001&0&114&88061\\
\bottomrule
\end{tabular}
\caption{Variables and constraints in the linear programs.}
\label{tab:lp-sizes}
\end{table}

\section{Coefficient data}\label{app:coefficient-data}

This appendix records the coefficient values used in the local inequalities.
Together with the definitions and verifications in
Appendix~\ref{app:checks}, these data suffice to check the local
inequalities and their averaged bounds.

The ordering conventions are fixed as follows. In the nine-point configuration, order
point labels lexicographically and take the lexicographically least label in
each vertex orbit, or the lexicographically least ordered pair in each
unordered pair orbit, as its representative. The table rows follow that
order. In the eighteen-point configuration, use the point order
$00,01,\ldots,08,10,11,\ldots,18$. The rows and columns of each coefficient
array are indexed by $0,1,\ldots,8$, and the squared-distance table is in
increasing order. For the square--tetrahedron configuration, the orbit
$\mathcal O_j$ is the orbit represented in row $j$ of the table below, and
that row order defines the indices used in $X_j$, $a_j$, and $B_q$.

\subsection{The $W_4$ inequality}\label{app:w4-pairs}

For the point order fixed in Section~\ref{sec:w4}, the pairs forming $\ell_2$ are
\begin{align*}
\mathcal E_4 = \{&(0,5),(0,15),(1,2),(1,4),(1,6),(1,16),\\
 &(2,3),(2,7),(2,17),(3,4),(3,8),(3,18),\\
 &(4,9),(4,19),(5,10),(6,7),(6,9),(6,11),\\
 &(7,8),(7,12),(8,9),(8,13),(9,14),(10,15),\\
 &(11,12),(11,14),(11,16),(12,13),(12,17),(13,14),\\
 &(13,18),(14,19),(16,17),(16,19),(17,18),(18,19)\}.
\end{align*}
The notation $(i,j)$ denotes the unordered pair $\{i,j\}$.

The dihedral symmetries of the unit square act on the twenty-point configuration in
the main text. To avoid repeated denominators, write
\[
10^6\mathcal L_{W_4}(r)
 = C+\sum_i A_i r_i+\sum_{i<j}B_{ij}r_ir_j.
\]
The constant is $C = -24994229$, and the tables below give the integers
$A_i$ and $B_{ij}$.
The linear coefficients are constant on the following three vertex orbits.

\begin{table}[H]
\centering
\scriptsize
\begin{tabular}{crr}
\toprule
indices & orbit size & $A_i$\\
\midrule
$0,5,10,15$ & $4$ & $3167237$\\
$1,2,7,8,13,14,16,19$ & $8$ & $5837166$\\
$3,4,6,9,11,12,17,18$ & $8$ & $10000000$\\
\bottomrule
\end{tabular}
\caption{Linear coefficients in the $W_4$ inequality.}
\label{w4:tab:linear}
\end{table}

For pairs, the first column gives the representative of a
pair orbit under these symmetries, the second gives its size, and the last
gives the common scaled integer $B_{ij}$. The corresponding coefficient in
$\mathcal L_{W_4}$ is $\beta_{ij} = 10^{-6}B_{ij}$. These rows determine
the complete upper-triangular matrix $(\beta_{ij})$.

\begin{table}[!htbp]
\centering
\scriptsize
\setlength{\tabcolsep}{5pt}
\begin{tabular}{crr}
\toprule
representative & orbit size & $B_{ij}$\\
\midrule
$(0,1)$ & $8$ & $3129119$\\
$(0,3)$ & $8$ & $-10000000$\\
$(0,5)$ & $4$ & $10000000$\\
$(0,6)$ & $8$ & $-2904617$\\
$(0,7)$ & $8$ & $355813$\\
$(0,8)$ & $8$ & $-2904617$\\
$(0,9)$ & $8$ & $-618433$\\
$(0,10)$ & $2$ & $10000000$\\
$(0,11)$ & $8$ & $-487123$\\
$(0,13)$ & $8$ & $9381567$\\
$(1,2)$ & $4$ & $10000000$\\
$(1,3)$ & $8$ & $-7095383$\\
$(1,4)$ & $8$ & $10000000$\\
$(1,6)$ & $8$ & $10000000$\\
$(1,7)$ & $8$ & $0$\\
$(1,8)$ & $4$ & $-8479162$\\
$(1,9)$ & $8$ & $-487123$\\
$(1,11)$ & $8$ & $-487123$\\
$(1,12)$ & $8$ & $-487123$\\
$(1,13)$ & $4$ & $4294145$\\
$(1,14)$ & $4$ & $4162834$\\
$(1,16)$ & $4$ & $10000000$\\
$(1,17)$ & $8$ & $0$\\
$(1,18)$ & $8$ & $2286183$\\
$(3,4)$ & $4$ & $10000000$\\
$(3,6)$ & $4$ & $355813$\\
$(3,9)$ & $8$ & $487123$\\
$(3,11)$ & $4$ & $0$\\
$(3,12)$ & $4$ & $131310$\\
$(3,18)$ & $4$ & $10000000$\\
\bottomrule
\end{tabular}
\caption{Pair coefficients in the $W_4$ inequality.}
\label{w4:tab:coefficients}
\end{table}

\begin{table}[!htbp]
\centering
\scriptsize
{\renewcommand{\arraystretch}{1.1}
\begin{tabular}{cr}
\toprule
$q$&$\sum_{\substack{i<j\\|v_i-v_j|^2 = q}}\beta_{ij}$\\
\midrule
$\tfrac12$&$-54.967048$\\
$1$&$360$\\
$\tfrac32+\sqrt2$&$-23.236936$\\
$\tfrac32$&$-2.10096$\\
$\tfrac32-\sqrt2$&$-23.236936$\\
$2$&$-40.660048$\\
$\tfrac52+\sqrt2$&$-3.896984$\\
$\tfrac52-\sqrt2$&$75.052536$\\
$2-\sqrt2$&$-3.896984$\\
$3-2\sqrt2$&$-17.265312$\\
$3$&$14.39248$\\
$4-2\sqrt2$&$17.17658$\\
$2+\sqrt2$&$3.896984$\\
$3+2\sqrt2$&$1.948492$\\
\bottomrule
\end{tabular}}
\caption{Nonzero coefficient sums for the $W_4$ kernel.}
\label{w4:tab:kernel}
\end{table}

With the kernel bound $K_{W_4}(t)-600J_0(t)<46$, the probability bound is
\[
\Phi_0(p) = -24.994229+185.366276p-342.794136p^2.
\]
The proof rounds these coefficients upwards to obtain $\Phi$.
Direct evaluation gives $\max_{0 \le p \le1/4}\Phi_0(p)<-0.07$.

\subsection{The nine-point inequality}

The group acting on the nine-point configuration has four vertex orbits and thirteen
orbits of unordered pairs. Each pair orbit is determined by the representative
in the table by applying the four symmetries of the configuration.

\begin{table}[H]
  \centering
  \begin{minipage}[t]{0.37\textwidth}
    \vspace{0pt}
    \centering
    \scriptsize
    \setlength{\tabcolsep}{3pt}
    \resizebox{\linewidth}{!}{%
    \begin{tabular}{ccc}
    \toprule
    representative & orbit size & coefficient $c_O$\\
    \midrule
    $00$ & $2$ & $-252$\\
    $01$ & $4$ & $-1004$\\
    $02$ & $2$ & $998$\\
    $11$ & $1$ & $-345$\\
    \bottomrule
    \end{tabular}%
    }
    \captionsetup{width = \linewidth,font = footnotesize,justification = raggedright,singlelinecheck = false}
    \captionof{table}{Nine-point vertex coefficients.}
    \label{tab:c3-vertex-coefficients}
  \end{minipage}\hfill
  \begin{minipage}[t]{0.61\textwidth}
    \vspace{0pt}
    \centering
    \scriptsize
    \setlength{\tabcolsep}{2pt}
    \begin{tabular}{cccrr}
    \toprule
    representative & orbit size & squared distance & $a_U$ & $b_U$\\
    \midrule
    $00,01$ & $4$ & $1$ & $-998$ & $2$\\
    $00,02$ & $4$ & $4$ & $-1000$ & $0$\\
    $00,11$ & $2$ & $3$ & $500$ & $150$\\
    $00,12$ & $4$ & $7$ & $1000$ & $300$\\
    $00,22$ & $1$ & $12$ & $375$ & $200$\\
    $01,02$ & $4$ & $1$ & $-446$ & $-894$\\
    $01,10$ & $2$ & $1$ & $2$ & $502$\\
    $01,11$ & $4$ & $1$ & $444$ & $948$\\
    $01,12$ & $2$ & $3$ & $662$ & $-12$\\
    $01,20$ & $4$ & $3$ & $448$ & $0$\\
    $01,21$ & $2$ & $4$ & $-264$ & $-12$\\
    $02,11$ & $2$ & $1$ & $-52$ & $-1000$\\
    $02,20$ & $1$ & $4$ & $112$ & $0$\\
    \bottomrule
    \end{tabular}
    \captionsetup{width = \linewidth,font = footnotesize}
    \captionof{table}{Nine-point pair coefficients.}
    \label{tab:c3-pair-coefficients}
  \end{minipage}
\end{table}

\subsection{The eighteen-point inequality}

The pair multiplicities and coefficient sums are recorded in
Table~\ref{l7:tab:distance-classes}. The next table gives the vertex
coefficients, followed by the three arrays containing all coefficients for
pairs.

The vertex coefficients are
\begin{equation*}
\begin{array}{c|rrrrrrrrr}
p&00&01&02&03&04&05&06&07&08\\
\hline
c_p&0&-54&-54&-16&-2&15&-46&-1&26\\[1mm]
p&10&11&12&13&14&15&16&17&18\\
\hline
c_p&18&-16&-52&1&-18&-23&-44&-53&0
\end{array}
\end{equation*}
For pairs, write
\[
\begin{aligned}
a^{00}_{jk}& = a_{(0,j),(0,k)},\\
a^{11}_{jk}& = a_{(1,j),(1,k)},\\
a^{01}_{jk}& = a_{(0,j),(1,k)}.
\end{aligned}
\]
The two within-column upper-triangular arrays and the cross-column array are
as follows. A dash denotes an unused diagonal or lower-triangular entry.
\begin{table}[H]
\centering
\begin{minipage}[t]{0.49\textwidth}
\centering
\scriptsize
\setlength{\tabcolsep}{1.5pt}
\begin{tabular}{crrrrrrrrr}
\toprule
$j\backslash k$&0&1&2&3&4&5&6&7&8\\
\midrule
0 & -- & $-100$ & $0$ & $0$ & $0$ & $0$ & $0$ & $0$ & $0$\\
1 & -- & -- & $-100$ & $-2$ & $-11$ & $-12$ & $14$ & $-16$ & $11$\\
2 & -- & -- & -- & $-100$ & $-54$ & $-100$ & $-58$ & $-22$ & $-43$\\
3 & -- & -- & -- & -- & $-100$ & $-95$ & $-100$ & $-60$ & $-55$\\
4 & -- & -- & -- & -- & -- & $-100$ & $-87$ & $-92$ & $-58$\\
5 & -- & -- & -- & -- & -- & -- & $-100$ & $-63$ & $-58$\\
6 & -- & -- & -- & -- & -- & -- & -- & $-100$ & $-42$\\
7 & -- & -- & -- & -- & -- & -- & -- & -- & $100$\\
8 & -- & -- & -- & -- & -- & -- & -- & -- & --\\
\bottomrule
\end{tabular}
\captionof{table}{The coefficients $a^{00}_{jk}$ for pairs in the first column.}
\end{minipage}\hfill
\begin{minipage}[t]{0.49\textwidth}
\centering
\scriptsize
\setlength{\tabcolsep}{1.5pt}
\begin{tabular}{crrrrrrrrr}
\toprule
$j\backslash k$&0&1&2&3&4&5&6&7&8\\
\midrule
0 & -- & $100$ & $-38$ & $-52$ & $-52$ & $-49$ & $-43$ & $5$ & $0$\\
1 & -- & -- & $100$ & $-57$ & $-77$ & $-46$ & $-13$ & $-12$ & $0$\\
2 & -- & -- & -- & $100$ & $-77$ & $-100$ & $-57$ & $20$ & $0$\\
3 & -- & -- & -- & -- & $100$ & $-95$ & $-100$ & $-12$ & $0$\\
4 & -- & -- & -- & -- & -- & $100$ & $-60$ & $-10$ & $0$\\
5 & -- & -- & -- & -- & -- & -- & $100$ & $1$ & $0$\\
6 & -- & -- & -- & -- & -- & -- & -- & $100$ & $0$\\
7 & -- & -- & -- & -- & -- & -- & -- & -- & $100$\\
8 & -- & -- & -- & -- & -- & -- & -- & -- & --\\
\bottomrule
\end{tabular}
\captionof{table}{The coefficients $a^{11}_{jk}$ for pairs in the second column.}
\end{minipage}
\end{table}
\begin{table}[H]
\centering
\scriptsize
\setlength{\tabcolsep}{3.2pt}
\begin{tabular}{crrrrrrrrr}
\toprule
$j\backslash k$&0&1&2&3&4&5&6&7&8\\
\midrule
0 & $-100$ & $0$ & $0$ & $0$ & $0$ & $0$ & $0$ & $0$ & $0$\\
1 & $-100$ & $-100$ & $40$ & $14$ & $31$ & $30$ & $57$ & $26$ & $0$\\
2 & $77$ & $-100$ & $-100$ & $99$ & $78$ & $79$ & $100$ & $59$ & $0$\\
3 & $63$ & $85$ & $-100$ & $-100$ & $93$ & $67$ & $68$ & $26$ & $0$\\
4 & $63$ & $57$ & $100$ & $-100$ & $-100$ & $84$ & $56$ & $28$ & $0$\\
5 & $52$ & $61$ & $65$ & $70$ & $-100$ & $-100$ & $85$ & $7$ & $0$\\
6 & $67$ & $96$ & $100$ & $73$ & $100$ & $-89$ & $100$ & $39$ & $0$\\
7 & $-85$ & $93$ & $90$ & $58$ & $60$ & $88$ & $100$ & $100$ & $0$\\
8 & $-90$ & $-80$ & $63$ & $47$ & $60$ & $60$ & $75$ & $100$ & $100$\\
\bottomrule
\end{tabular}
\caption{The coefficients $a^{01}_{jk}$ for pairs between the two columns.}
\end{table}
\subsection{The pentagon inequality}\label{app:p5-pairs}
Use the vectors $u_i$, translation $s$, and point labels $v_i$ from the
pentagon construction in Section~\ref{sec:p5}.
Within the two pentagons,
the side pairs are $(0,1),(1,2),(2,3),(3,4),(4,0)$ and
$(1,5),(5,6),(6,7),(7,8),(8,1)$.
For completeness, the calculation for pairs joining the two pentagons uses
\[
R_5^2 = \frac{5+\sqrt5}{10},\;
|u_i-(u_j+s)|^2 = |u_i-u_j-u_1+u_0|^2,
\]
and
\[
\cos 0 = 1,\;
\cos\frac{2\pi}{5} = \frac{\sqrt5-1}{4},\;
\cos\frac{4\pi}{5} = -\frac{\sqrt5+1}{4}.
\]
For $i\in\{0,1,2,3,4\}$ and $j\in\{1,2,3,4\}$, the pairs at distance one joining the
two pentagons are
$(i,j) = (1,1),(1,4),(2,2),(3,3),(4,4)$.
The first two are already the translated-pentagon sides $(1,5)$ and
$(1,8)$. The other three give $(2,6),(3,7),(4,8)$. 
Thus the complete
unit-pair list is
\begin{align*}
\mathcal E_5 = \{&(0,1),(0,4),(1,2),(1,5),(1,8),(2,3),(2,6),\\
 &(3,4),(3,7),(4,8),(5,6),(6,7),(7,8)\}.
\end{align*}
The point labels are ordered as in the pentagon configuration in the
main text. The linear coefficient vector is
\[
(\alpha_0,\ldots,\alpha_8) = (-23,77,-100,-23,-23,100,100,100,77).
\]
The coefficients for pairs are given by the following upper-triangular
matrix; blank entries are zero.
\begin{table}[H]
\centering
\scriptsize
\begin{tabular}{crrrrrrrrr}
\toprule
$i\backslash j$&0&1&2&3&4&5&6&7&8\\
\midrule
0&&$-100$&$23$&$23$&$-100$&&&&$23$\\
1&&&$-100$&$23$&$23$&$-100$&$-100$&$-100$&$-100$\\
2&&&&$-100$&$23$&$77$&$-100$&$77$&$100$\\
3&&&&&$-100$&&&$-100$&$23$\\
4&&&&&&&&&$-100$\\
5&&&&&&&$-100$&$-100$&$-100$\\
6&&&&&&&&$-48$&$-100$\\
7&&&&&&&&&$100$\\
8&&&&&&&&&\\
\bottomrule
\end{tabular}
\caption{Coefficients for pairs in the pentagon inequality.}
\label{p5:tab:coefficients}
\end{table}
\subsection{The four line inequalities}
The sequence labels in the following tables are binary. For $\phi_{1}$,
$\phi_{2}$, and $\phi_{3}$ they are the values of the red indicator
$r_j$. For $\phi_{4}$, they are the values
$\epsilon_j = r_j = (s_j+1)/2$. When the local weight is evaluated, the same
sequence is converted to signed values by $s_j = 2\epsilon_j-1$. Thus the tables
list only sequences satisfying the prescribed restrictions; these are the
entire domains of the functions $\phi_q$.
The two short tables of state values are
\begin{table}[H]
\centering
\begin{minipage}[t]{0.49\textwidth}
\centering
\scriptsize
\setlength{\tabcolsep}{2pt}
\begin{tabular}{cccccccc}
\toprule
sequence&$000$&$001$&$010$&$011$&$100$&$101$&$110$\\
\midrule
$\phi_{1}$&$0$&$-1$&$-2$&$-2$&$-3$&$-3$&$-3$\\
\bottomrule
\end{tabular}
\captionof{table}{Potential values for the $(\ell_3,\ell_4)$ inequality.}\label{hd:tab:phi34}
\end{minipage}\hfill
\begin{minipage}[t]{0.49\textwidth}
\centering
\scriptsize
\setlength{\tabcolsep}{2pt}
\begin{tabular}{crcrcrcr}
\toprule
seq.&$\phi_{4}$&seq.&$\phi_{4}$&seq.&$\phi_{4}$&seq.&$\phi_{4}$\\
\midrule
$0001$&$0$&$0010$&$-8$&$0011$&$-10$&$0100$&$-12$\\
$0101$&$-14$&$0110$&$-14$&$0111$&$-12$&$1000$&$-12$\\
$1001$&$-14$&$1010$&$-14$&$1011$&$-12$&$1100$&$-10$\\
$1101$&$-8$&$1110$&$0$&&&&\\
\bottomrule
\end{tabular}
\captionof{table}{Potential values for the $(\ell_4,\ell_4)$ inequality.}
\label{hd:tab:phi44}
\end{minipage}
\end{table}
The two longer tables of state values follow.
\begin{table}[H]
\centering
\begin{minipage}[t]{0.49\textwidth}
\centering
\scriptsize
\setlength{\tabcolsep}{2pt}
\begin{tabular}{@{}crcr@{}}
\toprule
sequence&$\phi_{2}$&sequence&$\phi_{2}$\\
\midrule
$000010$&$-2$&$000011$&$0$\\
$000100$&$-8$&$000101$&$-8$\\
$000110$&$-6$&$001000$&$-14$\\
$001001$&$-14$&$001010$&$-14$\\
$001011$&$-14$&$001100$&$-13$\\
$001101$&$-12$&$010000$&$-20$\\
$010001$&$-20$&$010010$&$-20$\\
$010011$&$-20$&$010100$&$-22$\\
$010101$&$-20$&$010110$&$-20$\\
$011000$&$-24$&$011001$&$-22$\\
$011010$&$-20$&$011011$&$-18$\\
$100001$&$-26$&$100010$&$-26$\\
$100011$&$-26$&$100100$&$-30$\\
$100101$&$-28$&$100110$&$-28$\\
$101000$&$-28$&$101001$&$-28$\\
$101010$&$-26$&$101011$&$-26$\\
$101100$&$-29$&$101101$&$-27$\\
$110000$&$-30$&$110001$&$-30$\\
$110010$&$-28$&$110011$&$-28$\\
$110100$&$-26$&$110101$&$-26$\\
$110110$&$-24$&&\\
\bottomrule
\end{tabular}
\captionof{table}{Potential values for the $(\ell_3,\ell_5)$ inequality.}
\label{hd:tab:phi35}
\end{minipage}\hfill
\begin{minipage}[t]{0.49\textwidth}
\centering
\scriptsize
\setlength{\tabcolsep}{2pt}
\begin{tabular}{@{}crcr@{}}
\toprule
sequence&$\phi_{3}$&sequence&$\phi_{3}$\\
\midrule
$000001$&$0$&$000010$&$-2$\\
$000011$&$-2$&$000100$&$-4$\\
$000101$&$-4$&$000110$&$-4$\\
$001000$&$-6$&$001001$&$-6$\\
$001010$&$-6$&$001011$&$-6$\\
$001100$&$-6$&$001101$&$-6$\\
$010000$&$-8$&$010001$&$-8$\\
$010010$&$-8$&$010011$&$-8$\\
$010100$&$-9$&$010101$&$-8$\\
$010110$&$-9$&$011000$&$-8$\\
$011001$&$-8$&$011010$&$-8$\\
$011011$&$-8$&$100000$&$-10$\\
$100001$&$-10$&$100010$&$-11$\\
$100011$&$-10$&$100100$&$-12$\\
$100101$&$-11$&$100110$&$-13$\\
$101000$&$-11$&$101001$&$-11$\\
$101010$&$-11$&$101011$&$-10$\\
$101100$&$-12$&$101101$&$-11$\\
$110000$&$-10$&$110001$&$-10$\\
$110010$&$-10$&$110011$&$-10$\\
$110100$&$-10$&$110101$&$-10$\\
$110110$&$-10$&&\\
\bottomrule
\end{tabular}
\captionof{table}{Potential values for the $(\ell_3,\ell_6)$ inequality.}
\label{hd:tab:phi36}
\end{minipage}
\end{table}
\subsection{The twenty-point square--tetrahedron inequality}\label{app:sq-labels}
The point order used throughout the coefficient tables is as follows.
Label the points by
\[
\begin{aligned}
&p_0 = 0,& &p_1 = w_0,& &p_2 = v_0,& &p_3 = v_0+w_0,\\
&p_4 = e,& &p_5 = e+w_0,& &p_6 = e+v_0,& &p_7 = e+v_0+w_0,\\
&p_8 = w_1,& &p_9 = v_1,& &p_{10} = v_1+w_1,& &p_{11} = e+w_1,\\
&p_{12} = e+v_1,& &p_{13} = e+v_1+w_1,& &p_{14} = w_2,& &p_{15} = v_2,\\
&p_{16} = v_2+w_2,& &p_{17} = e+w_2,& &p_{18} = e+v_2,&
&p_{19} = e+v_2+w_2.
\end{aligned}
\]
On the labels, the generators act as
\[
\begin{aligned}
\sigma = {}&(1\ 15)(2\ 14)(3\ 16)(5\ 18)(6\ 17)(7\ 19)
(8\ 9)(11\ 12),\\
\tau = {}&(0\ 4)(1\ 5)(2\ 6)(3\ 7)(8\ 11)(9\ 12)(10\ 13)
(14\ 17)(15\ 18)(16\ 19).
\end{aligned}
\]
In Table~\ref{sq:tab:pair-orbits}, the representative is a pair of point labels from the
square--tetrahedron section. The third column is the size of its orbit under
the group $G$, the fourth is the common squared distance, and the last is the
coefficient $a_j$ in the square--tetrahedron expression.
\begin{table}[H]
\centering
\scriptsize
\begin{tabular}{cr cr}
\toprule
$q$&$B_q$&$q$&$B_q$\\
\midrule
$2-\sqrt3$&$-68300$&$1$&$-300708$\\
$3-\sqrt3$&$56588$&$2$&$-142872$\\
$4-\sqrt3$&$-14632$&$3$&$18334$\\
$5-\sqrt3$&$10996$&$2+\sqrt3$&$628$\\
$4$&$-5700$&$3+\sqrt3$&$-404$\\
$4+\sqrt3$&$204$&$6$&$386$\\
$5+\sqrt3$&$380$&$7$&$136$\\
\bottomrule
\end{tabular}
\caption{Coefficient sums at each squared distance in the square--tetrahedron
kernel.}
\label{sq:tab:distance-coefficients}
\end{table}
\begin{table}[H]
\centering
\scriptsize
\setlength{\tabcolsep}{2.5pt}
\caption{Pair orbits and integer coefficients for the square--tetrahedron inequality.}
\label{sq:tab:pair-orbits}
\begin{minipage}[t]{0.47\textwidth}
\centering
\begin{tabular}{@{}r c r c r@{}}
\toprule
$j$&representative&$|\mathcal O_j|$&$q_j$&$a_j$\\
\midrule
0&$\{8,10\}$&4&$1$&$0$\\
1&$\{0,8\}$&4&$1$&$-17344$\\
2&$\{0,10\}$&2&$2$&$-152$\\
3&$\{0,1\}$&4&$1$&$-9839$\\
4&$\{4,8\}$&4&$2$&$-21087$\\
5&$\{0,2\}$&4&$1$&$-10680$\\
6&$\{4,10\}$&2&$3$&$587$\\
7&$\{0,3\}$&4&$2$&$-2319$\\
8&$\{3,4\}$&4&$3$&$2040$\\
9&$\{1,4\}$&4&$2$&$-553$\\
10&$\{1,3\}$&4&$1$&$3317$\\
11&$\{2,3\}$&4&$1$&$-997$\\
12&$\{7,9\}$&4&$3-\sqrt3$&$5100$\\
13&$\{0,4\}$&1&$1$&$-39099$\\
14&$\{5,14\}$&4&$4$&$-1425$\\
15&$\{2,4\}$&4&$2$&$-1277$\\
16&$\{2,14\}$&2&$2+\sqrt3$&$-12$\\
17&$\{1,14\}$&4&$3$&$1304$\\
18&$\{3,14\}$&4&$4+\sqrt3$&$51$\\
19&$\{6,14\}$&2&$3+\sqrt3$&$-22$\\
20&$\{1,9\}$&4&$2-\sqrt3$&$-6903$\\
21&$\{7,14\}$&4&$5+\sqrt3$&$95$\\
22&$\{1,2\}$&4&$2$&$-1171$\\
23&$\{5,9\}$&4&$3-\sqrt3$&$5627$\\
24&$\{1,5\}$&2&$1$&$-6775$\\
25&$\{5,10\}$&4&$3-\sqrt3$&$485$\\
26&$\{2,5\}$&4&$3$&$1144$\\
27&$\{2,9\}$&4&$1$&$-9068$\\
\bottomrule
\end{tabular}
\end{minipage}%
\hspace{1.5em}%
\begin{minipage}[t]{0.47\textwidth}
\centering
\begin{tabular}{@{}r c r c r@{}}
\toprule
$j$&representative&$|\mathcal O_j|$&$q_j$&$a_j$\\
\midrule
28&$\{1,8\}$&4&$1$&$-2603$\\
29&$\{3,5\}$&4&$2$&$-2761$\\
30&$\{3,15\}$&4&$4-\sqrt3$&$-3658$\\
31&$\{1,10\}$&4&$2-\sqrt3$&$-337$\\
32&$\{3,16\}$&2&$6$&$193$\\
33&$\{7,15\}$&4&$5-\sqrt3$&$2749$\\
34&$\{7,16\}$&2&$7$&$68$\\
35&$\{2,10\}$&4&$2+\sqrt3$&$-294$\\
36&$\{3,9\}$&4&$2-\sqrt3$&$-6010$\\
37&$\{3,7\}$&2&$1$&$-2480$\\
38&$\{9,11\}$&2&$3$&$460$\\
39&$\{5,8\}$&4&$2$&$-5565$\\
40&$\{1,15\}$&2&$2-\sqrt3$&$-7650$\\
41&$\{6,8\}$&4&$3+\sqrt3$&$-21$\\
42&$\{7,8\}$&4&$3+\sqrt3$&$-145$\\
43&$\{10,11\}$&4&$2$&$537$\\
44&$\{2,8\}$&4&$2+\sqrt3$&$-93$\\
45&$\{3,10\}$&4&$2$&$-59$\\
46&$\{3,6\}$&4&$2$&$1433$\\
47&$\{6,9\}$&4&$2$&$-2065$\\
48&$\{3,8\}$&4&$2+\sqrt3$&$550$\\
49&$\{5,15\}$&2&$3-\sqrt3$&$5870$\\
50&$\{6,10\}$&4&$3+\sqrt3$&$76$\\
51&$\{8,9\}$&2&$2$&$-1510$\\
52&$\{7,10\}$&4&$3$&$-428$\\
53&$\{10,13\}$&1&$1$&$-267$\\
54&$\{2,6\}$&2&$1$&$-54$\\
55&$\{8,11\}$&2&$1$&$-26934$\\
\bottomrule
\end{tabular}
\end{minipage}
\end{table}


\begin{thebibliography}{99}
\bibitem{Ambrus2024}
G. Ambrus, A. Csisz\'arik, M. Matolcsi, D. Varga, and P. Zs\'amboki,
\emph{The density of planar sets avoiding unit distances},
Math. Program. \textbf{207} (2024), no. 1--2, 303--327.
\bibitem{Arman2018}
A. Arman and S. Tsaturian,
\emph{A result in asymmetric Euclidean Ramsey theory},
Discrete Math. \textbf{341} (2018), no. 5, 1502--1508.
\bibitem{Bekka2008}
B. Bekka, P. de la Harpe, and A. Valette,
\emph{Kazhdan's Property~(T)}, New Mathematical Monographs, vol.~11,
Cambridge University Press, Cambridge, 2008.
\bibitem{Billingsley1999}
P. Billingsley, \emph{Convergence of Probability Measures}, 2nd ed.,
Wiley, New York, 1999.
\bibitem{Cantwell1996}
K. Cantwell,
\emph{Finite Euclidean Ramsey theory},
J. Combin. Theory Ser. A \textbf{73} (1996), no. 2, 273--285.
\bibitem{Conlon2019}
D. Conlon and J. Fox,
\emph{Lines in Euclidean Ramsey theory},
Discrete Comput. Geom. \textbf{61} (2019), no. 1, 218--225.
\bibitem{ConlonFuhrer2026}
D. Conlon and J. F\"uhrer,
\emph{Non-spherical sets versus lines in Euclidean Ramsey theory},
Canad. Math. Bull. \textbf{69} (2026), no. 1, 179--183.
\url{https://doi.org/10.4153/S0008439525101082}.
\bibitem{ConlonWu}
D. Conlon and Y.-H. Wu,
\emph{More on lines in Euclidean Ramsey theory},
C. R. Math. Acad. Sci. Paris \textbf{361} (2023), 897--901.
\bibitem{Csizmadia1994}
G. Csizmadia and G. T\'oth,
\emph{Note on a Ramsey-type problem in geometry},
J. Combin. Theory Ser. A \textbf{65} (1994), no. 2, 302--306.
\bibitem{CurrierModyXieZhang2026}
G. Currier, P. Mody, Z. Xie, and J. Zhang,
\emph{Improved bounds for lines and $1$-separated sets in Euclidean Ramsey
theory}, arXiv:2606.17194 [math.CO], 2026.
\bibitem{CurrierMooreYip3AP}
G. Currier, K. Moore, and C. H. Yip,
\emph{Any two-coloring of the plane contains monochromatic 3-term arithmetic progressions},
Combinatorica \textbf{44} (2024), no. 6, 1367--1380.
\bibitem{CurrierMooreYipAvoiding}
G. Currier, K. Moore, and C. H. Yip,
\emph{Avoiding short progressions in Euclidean Ramsey theory},
J. Combin. Theory Ser. A \textbf{217} (2026), 106080.
\bibitem{deGrey2018}
A. D. N. J. de Grey,
\emph{The chromatic number of the plane is at least 5},
Geombinatorics \textbf{28} (2018), no. 1, 18--31.
\bibitem{deOliveira2010}
F.~M. de Oliveira Filho and F. Vallentin,
\emph{Fourier analysis, linear programming, and densities of distance avoiding sets in $\mathbb{R}^{n}$},
J. Eur. Math. Soc. \textbf{12} (2010), no. 6, 1417--1428.
\bibitem{DiaconisFreedman1987}
P.~Diaconis and D.~Freedman,
\emph{A Dozen de Finetti-Style Results in Search of a Theory},
Ann. Inst. Henri Poincar\'e Probab. Stat. \textbf{23} (1987), no. S2, 397--423.
\bibitem{Erdos1973}
P. Erd\H{o}s, R. L. Graham, P. Montgomery, B. L. Rothschild,
J. Spencer, and E. G. Straus,
\emph{Euclidean Ramsey theorems. I},
J. Combin. Theory Ser. A \textbf{14} (1973), no. 3, 341--363.
\bibitem{Erdos1975a}
P. Erd\H{o}s, R. L. Graham, P. Montgomery, B. L. Rothschild,
J. H. Spencer, and E. G. Straus,
\emph{Euclidean Ramsey theorems II},
in \emph{Infinite and Finite Sets},
Colloquia Mathematica Societatis J\'anos Bolyai, vol. 10,
North-Holland, Amsterdam, 1975, pp. 529--557.
\bibitem{Erdos1975b}
P. Erd\H{o}s, R. L. Graham, P. Montgomery, B. L. Rothschild,
J. H. Spencer, and E. G. Straus,
\emph{Euclidean Ramsey theorems III},
in \emph{Infinite and Finite Sets},
Colloquia Mathematica Societatis J\'anos Bolyai, vol. 10,
North-Holland, Amsterdam, 1975, pp. 559--583.
\bibitem{ErdosGraham1979}
P. Erd\H{o}s and R. L. Graham,
\emph{Old and new problems and results in combinatorial number theory:
van der Waerden's theorem and related topics},
Enseign. Math. (2) \textbf{25} (1979), no. 3--4, 325--344.
\bibitem{Falconer1981}
K.~J. Falconer,
\emph{The realization of distances in measurable subsets covering $\mathbb R^n$},
J. Combin. Theory Ser. A \textbf{31} (1981), no.~2, 184--189.
\bibitem{Folland2016}
G.~B. Folland,
\emph{A Course in Abstract Harmonic Analysis}, 2nd ed.,
Textbooks in Mathematics, CRC Press, Boca Raton, FL, 2016.
\bibitem{FuhrerToth}
J. F\"uhrer and G. T\'oth,
\emph{Progressions in Euclidean Ramsey theory},
European J. Combin. \textbf{125} (2025), 104105.
\bibitem{Gasarch2025}
W. Gasarch, A. Gezalyan, and R. Parker,
\emph{Open Problems Column: Monochromatic Unit Squares: Exposition and Open Problems},
SIGACT News \textbf{56} (2025), no. 3, 38--55,
\url{https://doi.org/10.1145/3767145.3767149}.
\bibitem{HuLiuToAppear}
G. Hu and Y. Liu,
\emph{A Two-Coloring of the Plane Avoiding a Red $\ell_2$ and Blue $\ell_{42}$},
to appear.
\bibitem{Ivan1979}
L. Iv\'an,
\emph{Monochromatic point sets in the plane and in the space},
Master's thesis, University of Szeged, Bolyai Institute, 1979 (in Hungarian).
\bibitem{Juhasz1979}
R. Juh\'asz,
\emph{Ramsey type theorems in the plane},
J. Combin. Theory Ser. A \textbf{27} (1979), no. 2, 152--160.
\bibitem{Moser1961}
L. Moser and W. Moser,
\emph{Solution to problem 10},
Canad. Math. Bull. \textbf{4} (1961), no. 2, 187--189.
\bibitem{DLMF2026}
F.~W.~J. Olver, A.~B. Olde Daalhuis, D.~W. Lozier,
B.~I. Schneider, R.~F. Boisvert, C.~W. Clark,
B.~R. Miller, B.~V. Saunders, H.~S. Cohl, and M.~A. McClain (eds.),
\emph{NIST Digital Library of Mathematical Functions},
\url{https://dlmf.nist.gov/}, Release 1.2.8 of 2026-09-15.
\bibitem{OpenAI}
OpenAI,
\emph{The Euclidean plane is not five-colorable},
preprint, September 23, 2026.
\href{https://github.com/openai/math/blob/fd4aeeb2ee4fc729c18d98444fed42fd0529eeeb/preprints/The-Euclidean-plane-is-not-five-colorable-September-23-2026/paper.pdf}{OpenAI mathematics repository}.
\bibitem{Payne2009}
M.~S. Payne,
\emph{Unit distance graphs with ambiguous chromatic number},
Electron. J. Combin. \textbf{16} (2009), no.~1, N31.
\bibitem{Szekely1984}
L.~A. Sz\'ekely,
\emph{Measurable chromatic number of geometric graphs and sets without
some distances in Euclidean space},
Combinatorica \textbf{4} (1984), no. 2--3, 213--218.
\bibitem{Szlam2001}
A. D. Szlam,
\emph{Monochromatic translates of configurations in the plane},
J. Combin. Theory Ser. A \textbf{93} (2001), no. 1, 173--176.
\bibitem{Tsaturian2017}
S. Tsaturian,
\emph{A Euclidean Ramsey result in the plane},
Electron. J. Combin. \textbf{24} (2017), no. 4, Paper No. P4.35, 9 pp.
\end{thebibliography}
\end{document}